%% file: CMFG.tex
\documentclass[11pt]{article}

\usepackage{fullpage}
\usepackage[utf8]{inputenc}
\usepackage{algorithm, algpseudocode}
\usepackage[T1]{fontenc}  
\usepackage{todonotes}
\usepackage{amsmath}
\usepackage{amsfonts,amssymb}
\usepackage[margin=0.75in]{geometry}
\usepackage{url,hyperref}
\usepackage{amsthm}
\usepackage{xcolor}
\usepackage{float}
\usepackage{booktabs} 
\usepackage{graphicx}
\usepackage{caption}
\usepackage{subcaption}
\usepackage{ulem}
\usepackage{bm}
\usepackage[shortlabels]{enumitem}

\allowdisplaybreaks

\input defs_cmfg.tex

\newcommand{\norm}[1]{{\left\vert\kern-0.25ex\left\vert\kern-0.25ex\left\vert #1
    \right\vert\kern-0.25ex\right\vert\kern-0.25ex\right\vert}}

\author{
Anran Hu 
\thanks{Columbia University, IEOR.
\textbf{Email:} \texttt{ah4277@columbia.edu}}
\and
Zijiu Lyu
\thanks{UC Berkeley, IEOR. \textbf{Email:} \texttt{zijiu.lyu@berkeley.edu}} 
}
\title{Mean-Field Games with Constraints}
\date{}
\begin{document}

\maketitle

\begin{abstract}
This paper introduces a framework of Constrained Mean-Field Games (CMFGs), where each agent solves a constrained Markov decision process (CMDP). This formulation captures scenarios in which agents' strategies are subject to feasibility, safety, or regulatory restrictions, thereby extending the scope of classical mean field game (MFG) models. We first establish the existence of CMFG equilibria under a strict feasibility assumption, and we further show uniqueness under a classical monotonicity condition. To compute equilibria, we develop Constrained Mean-Field Occupation Measure Optimization (CMFOMO), an optimization-based scheme that parameterizes occupation measures and shows that finding CMFG equilibria is equivalent to solving a single optimization problem with convex constraints and bounded variables. CMFOMO does not rely on uniqueness of the equilibria and can approximate all equilibria with arbitrary accuracy. We further prove that CMFG equilibria induce $O(1 \slash \sqrt{N})$-Nash equilibria in the associated constrained $N$-player games, thereby extending the classical justification of MFGs as approximations for large but finite systems. Numerical experiments on a modified Susceptible-Infected-Susceptible (SIS) epidemic model with various constraints illustrate the effectiveness and flexibility of the framework.
\end{abstract}

% \textcolor{red}{
% % Urich horst’s papers on control/game with trading constraints (and references therein).
% % Stockbridge LP for control problems with constraints.
% equality constraints?
% }

\section{Introduction}
In the domain of game theory, the analysis of large-scale strategic interactions poses both profound opportunities and formidable challenges. Multi-agent systems, in which numerous agents interact under shared decision-making rules, are central to fields such as economics, finance, and engineering \cite{yang2020overview,zhang2021multi}. Yet the complexity of these systems grows rapidly with the number of participants, making the direct study of $N$-player stochastic games analytically and computationally intractable.

Introduced by the seminal works of \cite{huang2006large} and \cite{lasry2007mean}, mean field games (MFGs) provide a powerful framework to address this difficulty. The key idea is to approximate large-population games by considering the limiting regime as the number of agents tends to infinity. In this regime, the influence of any single agent becomes negligible, and each agent interacts only with the collective effect of the population. This simplification transforms an otherwise intractable many-player game into a tractable model where equilibrium analysis can be carried out systematically.

MFG theory has not only offered conceptual clarity but also rigorous approximation results: equilibria of the limiting mean field game can be used to construct approximate Nash equilibria for the corresponding $N$-player game, with explicit error bounds of order $O(1/\sqrt{N})$ \cite{huang2006large,lasry2007mean,delarue2019master,saldi2018markov}. Over the past decade, this approach has been extensively studied and extended, leading to a rich literature on existence, uniqueness, and numerical methods. Its versatility has enabled applications across domains ranging from financial markets and systemic risk to crowd dynamics, epidemiology, and large-scale social networks.

Despite these successes, classical MFG models typically assume that agents operate without binding restrictions on their states or controls. This unconstrained setting is analytically tractable but often unrealistic. In practice, agents’ decisions are almost always shaped by feasibility, safety, or regulatory requirements. The relevance of constraints has long been recognized in related areas such as constrained Markov decision processes (CMDPs), which capture situations where an individual agent face multiple competing objectives. Examples include telecommunication, where throughput must be maximized while keeping delays below a threshold \cite{altmanCMDP}; finance, where risk and return trade-offs are enforced through constraints as in Markowitz’s portfolio theory \cite{Markowitz}; and resource allocation problems in medical systems, where limited capacity dictates feasible allocations \cite{medical}. More recently, the study of safe reinforcement learning has emphasized constraints as a safeguard against catastrophic outcomes \cite{wachi2020safe}.

While classical MFG models have yielded many elegant results, they often overlook constraints that naturally arise in real-world systems. Incorporating such constraints can make the framework more realistic and broaden its range of applications.  At the same time, constraints introduce significant analytical challenges: the feasible set of agent decisions may depend on the population distribution itself, which prevents a direct application of existing methods and analyses.

\paragraph{Our work and approach.} In this paper, we introduce a framework of Constrained Mean Field Games (CMFGs), where both the time horizon and the state space are discrete and finite. Unlike the classical MFG setting, where the representative agent solves an unconstrained Markov decision process (MDP), here each agent faces a constrained MDP (CMDP) that maximizes utility subject to restrictions on states and actions. In the general case, constraints depend on an individual’s state and action together with the population distribution, a setting we call agent-population-level constraints. We define CMFG Nash equilibria through a fixed-point argument: given a population distribution, the representative agent solves the induced CMDP to obtain its best response, and this policy must in turn regenerate the same distribution. For clarity, we use the term MFG to refer exclusively to the unconstrained case, and the terms player and agent interchangeably.

For CMFGs with general agent-population-level constraints, we prove the existence of equilibria under suitable assumptions using Kakutani’s fixed-point theorem, in line with the classical MFG literature. A key difference is that, with constraints, the feasible domain of each agent depends on the population distribution, which complicates the analysis. To apply Kakutani’s theorem, we impose a strict feasibility assumption: for every constraint and any mean-field distribution, there must exist a policy that respects all constraints while satisfying the chosen one with a positive margin. This guarantees that the feasible set for each agent is nonempty and has enough slack to avoid degenerate boundary cases. In addition to existence, we also provide a uniqueness result under stronger assumptions. Specifically, when the transition probabilities and constraints are independent of the mean-field distribution and the reward functions satisfy the classical monotonicity condition of Lasry and Lions, the CMFG admits a unique equilibrium. This complements our general existence result and parallels the uniqueness theory developed for classical MFGs.

To compute equilibria under agent-population-level constraints, we develop a numerical scheme called Constrained Mean Field Occupation Measure Optimization (CMFOMO). The method parameterizes the agent’s occupation measure and leverages Karush-Kuhn-Tucker (KKT) conditions to reformulate the equilibrium problem as an optimization program with convex constraints and bounded variables. CMFOMO does not assume uniqueness and can characterize all CMFG equilibria. We show that suboptimal solutions of CMFOMO approximate true equilibria arbitrarily well when the objective is minimized sufficiently.

We then consider a population-level constraint setting, a degenerate case where constraints depend only on the mean-field distribution and not on individual strategies, without imposing the strict feasibility assumption. In this case, the problem reduces to finding standard Nash equilibria that satisfy the population-level constraints, which can be solved using a modified CMFOMO framework by removing the components associated with the strict feasibility assumption, and similar approximation guarantees are obtained.

Finally, we establish a rigorous connection between CMFGs and finite-player constrained games. In particular, we show that CMFG equilibria yield $O(1/\sqrt{N})$-Nash equilibria for the corresponding constrained $N$-player games. This extends the classical justification of MFGs as approximations for large but finite-player games to the constrained setting. 

In the end, numerical experiments on a modified Susceptible–Infected–Susceptible (SIS) model with several different kinds of constraints illustrate the effectiveness of our approach in both the agent-population-level and population-level settings.

\paragraph{Related work.} Our approach builds on the linear programming (LP) method for CMDPs, where decision variables are occupation measures. We refer to \cite{altmanCMDP} for a comprehensive treatment of solution approaches and theory, and to \cite{altman1991markov,altman1994denumerable,altmanCMDP3rd} for early developments. %Another closely related concept is Safe Reinforcement Learning (Safe RL), where the agent not only maximizes its reward but also keeps track of its current status to avoid breaking the so-called safety constraints. And it has been a common practice to attack such problems by the tools of CMDPs. We refer to \cite{wachi2020safe, thomas2021safe, wagener2021safe, bura2022dope, sootla2022saute} for some recent papers on this topic. 
In the continuous-time setting, the LP method has also been applied to stochastic control problems, potentially with constraints. See \cite{kurtz2017linear} and the references therein for a broad survey.

In recent years, a growing body of work has used LP formulations and occupation measures in games without agent constraints. For example, \cite{bouveret2020mean,dumitrescu2021control} introduce LP formulations for continuous-time MFGs to characterize optimality conditions and establish existence of Nash equilibria with mixed strategies. %Building on this perspective, 
\cite{mfomo} develops a primal-dual characterization for discrete-time MFGs, leading to an optimization-based formulation (MFOMO), which is later extended to continuous-time MFGs in \cite{guo2025continuous}. Using occupation measures, \cite{hu2024mf} proposes a globally convergent algorithm for computing equilibria of MFGs under the monotonicity condition, and further extended the analysis to reinforcement learning with regret guarantees in multi-agent systems. Our CMFOMO generalizes the MFOMO framework of \cite{mfomo} by incorporating agent constraints while preserving its favorable properties.

% Imposing agent constraints in MFGs has also been studied from different perspectives. For instance, \cite{cannarsa2018existence} considers continuous-time deterministic MFGs where the agent state is confined into a compact subset of the Euclidean space (which are called the state constraints), and proves the existence and uniqueness of the NE in the corresponding context. This problem is then studied in \cite{cannarsa2021mean, cannarsa2021weak, arjmand2022nonsmooth} using PDE approach and the results are extended by \cite{graber2021note,bonnans2023lagrangian} to incorporate mean-field game of controls and mixed state-action constraints. %Besides, \cite{bonnans2023lagrangian} studies a deterministic mean-field game of control with both state constraints and mixed state-action constraints.
% In the stochastic games, \cite{fu2021mean} studies portfolio liquidation MFGs where terminal inventory has to be zero \cite{fu2021mean}, where \cite{fu2024mean} focuses the constraint where players can only trade in a single direction. For the finite player setting, \cite{evangelista2020finite} studies a liquidation game with terminal liquidation constraint and establishes its relation with the corresponding MFG. Compared to the existing works, this paper investigate stochastic MFGs with more general constraints, which can depend on single player's states, actions and the population distributions. We also provide an optimization scheme to solve the constrained MFGs, and analyze the approximation error between MFGs and $N$-player games.

Imposing agent constraints in MFGs has been investigated from several perspectives. In deterministic settings, \cite{cannarsa2018existence} considers continuous-time MFGs where the agent’s state is confined to a compact subset of the Euclidean space, which is referred to as state constraints, and establishes existence and uniqueness of Nash equilibria. Subsequent works such as \cite{cannarsa2021mean, cannarsa2021weak, arjmand2022nonsmooth} extend this analysis using PDE-based methods, while \cite{graber2021note, bonnans2023lagrangian} further generalize the framework to include mean-field games of controls and mixed state–action constraints.
In stochastic formulations, \cite{fu2021mean} studies portfolio liquidation MFGs with a terminal inventory constraint, requiring each agent’s position to be fully liquidated by the terminal time. Building on this, \cite{fu2024mean} introduces directional trading constraints, where each player is restricted to trade in a single direction depending on their initial position. For finite-player settings, \cite{evangelista2020finite} analyzes a stochastic liquidation game with terminal constraints and rigorously connects its equilibria to those of the corresponding MFG.

Compared to these studies, our paper investigates stochastic MFGs with more general constraints, which may depend simultaneously on individual agents’ states, actions, and the population distribution. In addition, we propose an optimization-based computational framework for solving such constrained MFGs and provide theoretical guarantees on the approximation error between MFG and $N$-player equilibria.

%Constraints also emerge in the context of portfolio liquidation games, \cite{fu2021mean,fu2020mean} study MFGs and stackelberg MFGs with terminal state constraint

% portfolio liquidation games:
% \cite{evangelista2020finite} finite player games, existence and uniqueness, characterization

% (choose a subset of the followings)
% \cite{fu2021mean} terminal state constraint, MFG
% %\cite{fu2020mean} terminal state constraint, stackelberg MFG
% \cite{fu2024mean} trade in a single direction
% %\cite{fu2022portfolio} Hawkes
% \cite{fu2024mean2} market drop-out, absorption constraint

% other games:
% % \cite{bonnans2023lagrangian} state constraint and mixed state-action constraint, deterministic mean-field game of control
% \cite{graewe2022maximum} maximum principle, deterministic mfg, absorption constraint

% % Besides state constraints, \cite{fu2024mean} considers the MFG of optimal portfolio liquidation where the players (\text{i.e.}, brokers) are restricted to trade in a single direction (\text{i.e.}, either long or short), depending on their initial positions. \cite{fu2020mean, fu2024mean2} are also on similar topics.

% ======

% To clarify, their settings are essentially different from ours (though sharing a similar name). The proof of the NE uniqueness in \cite{cannarsa2018existence} relies on the assumption that the feasible domain of agent state is unaffected by the population distribution, which we do not assume in our setup. For another example, 

% XXX. 

\paragraph{Organization of the paper.} The remainder of the paper is organized as follows. Section~\ref{sec: setup} presents the formulation of CMFGs and introduces both agent-population-level and population-level constraints. Section~\ref{sec: agent-population-level} develops existence and uniqueness results as well as our numerical scheme, CMFOMO. Section~\ref{sec: pop-level} treats the population-level setting and its approximation results. Section~\ref{sec: N_player} establishes the connection to constrained $N$-player games, and Section~\ref{sec: numerical} illustrates the framework with numerical experiments on the SIS model.

\section{Constrained Mean-Field Games} \label{sec: setup}

In this section, we formally introduce (Markov) Constrained Mean-Field Games (CMFGs). The game is defined over a finite time horizon $\mathcal{T} := \{0, 1, \dots, T\}$, where $T < +\infty$, with a finite state space $\mathcal{S} \subset \mathbb{R}^d$ and a finite action space $\mathcal{A} \subset \mathbb{R}^n$. In this game, there is an infinite population of homogeneous players, each independently starting from an initial state $s_0 \sim \mu_0 \in \Delta(\mathcal{S})$. At each time step $t \in \mathcal{T}$, a representative player observes their current state $s_t$ and selects an action $a_t \in \mathcal{A}$ according to a randomized (Markov) policy $\pi_t(\cdot \,|\, s_t) \in \Delta(\mathcal{A})$. The player then receives a reward $r_t(s_t, a_t, L_t)$ and transitions to a new state $s_{t+1}$ following the transition probability $p_t(\cdot \,|\, s_t, a_t, L_t)$, where $L_t \in \Delta(\mathcal{S} \times \mathcal{A})$ represents the joint state-action distribution of the population at time $t$.

This formulation aligns with the standard framework of discrete time Mean-Field Games (MFGs), where the players are free to optimize their policies without additional restrictions. However, in CMFGs, players are subject to certain policy constraints. Alongside receiving rewards, players also incur a cost $c_t(s_t, a_t, L_t)$, defined by a cost function $c: \mathcal{T} \times \mathcal{S} \times \mathcal{A} \times \Delta(\mathcal{S} \times \mathcal{A}) \to \mathbb{R}^k$. Each player must ensure that their expected cumulative cost over the time horizon does not exceed a specified threshold $\gamma_0 \in \mathbb{R}^k$. This constraint enforces that the expected cumulative cost remains within the prescribed limits, introducing an additional layer of complexity to the decision-making process.

Mathematically, given a mean-field flow $L := {L_t \in \Delta(\state \times \action)}_{t \in \T}$ that characterizes the joint distribution of states and actions in the population at each time step, each player solves the following Constrained Markov Decision Process (CMDP). % Importantly, as the number of players approaches infinity, they become asymptotically independent\footnote{This independence arises due to the mean-field approximation, which decouples individual dynamics from direct interactions with other players.}.

\begin{cmdps} \label{problem: cmdps}
    \begin{align}
        \text{maximize}_{\{ \pi_t \}_{t \in \T}} & \quad \E \left[ \sum_{t \in \T} r_t(s_t, a_t, L_t) \,\bigg|\, s_0 \sim \mu_0 \right] \label{eq: cmdp1} \tag{CMDP reward} \\
        \text{subject to} & \quad a_t \sim \pi_t(s_t), \quad s_{t + 1} \sim p_t(s_t, a_t, L_t), \tag{CMDP dynamics} \\
        & \quad \E \left[ \sum_{t \in \T} c_t(s_t, a_t, L_t) \,\bigg|\, s_0 \sim \mu_0 \right] \leq \gamma_0, \label{eq: cmdp3} \tag{CMDP constraints}
    \end{align}
\end{cmdps}
\noindent
where $\mu_0 \in \Delta(\state)$ denotes the initial distribution. For notational convenience, let $\cmdp{L}$ denote the constrained MDP problem described above, and let $\optimal{L}$ denote its optimal expected cumulative reward, whenever it exists.

The extra constraints in \eqref{eq: cmdp3} show the main difference of our CMFGs compared with the classic unconstrained MFGs and one needs to define a new solution concept that is tailored to our setting. To formally define the Nash Equilibrium (NE) of CMFGs, we first define the mapping $\Psi(\cdot) \in \big[\Delta(\state \times \action)\big]^{|\T|}$ which maps any admissible policy $\pi$ to the flow of population state-action distribution induced by $\pi$. The mapping $\Psi$ is defined recursively in the following way. For any $s \in \state$ and $a \in \action$, 
%we denote by $\Psi(\pi) \in [\Delta(\state \times \action)]^{|\T|}$ the flow of population distribution induced solely by $\pi$, \textit{i.e.}, for each $s \in \state$ and $a \in \action$, we have:
\begin{align*}
    \begin{cases}
         \Psi_0(\pi) (s, a) := \pi_0(a \,|\, s) \mu_0(s), \\
        \Psi_{t + 1}(\pi) (s, a) := \pi_{t+1}(a\,|\,s)\sum_{s'\in\mathcal{S},a'\in\mathcal{A}}p_t(s \,|\, s',a',\Psi_{t}(\pi))\Psi_{t}(\pi)(s',a'). 
        % \Big[\Psi_{t}(\pi) \cdot p_t(s\,|\, \Psi_{t}(\pi))\Big] \pi_{t + 1}(s)(a).
    \end{cases}
\end{align*}

\begin{definition}[Nash equilibrium of CMFGs] \label{def: constrained ne}
    We say that $(\pi, L)$ is a Nash equilibrium of CMFG if it satisfies:
    \begin{enumerate}[1)]
        \item $\pi$ is an optimal policy of the constrained MDP $\cmdp{L}$;
        \item $L = \Psi(\pi)$.
    \end{enumerate}
\end{definition}

\begin{remark} \label{rmk: opt_consis}
    In the literature of MFGs, conditions 1) and 2) are often referred to as the \textit{optimality/best response} and \textit{consistency/compatibility} conditions, respectively. Here,  the condition 1) involves a constrained MDP problem instead of a standard MDP problem, which is the key difference of our framework compared with the standard definition of the NEs for MFGs. The NE of CMFGs indicates that no agent can  improve its payoff by unilaterally deviating from current strategy while satisfying constraints.
\end{remark}

\paragraph{Population-level constraints and agent-population-level constraints.} 
When the cost functions take different forms, two types of constraints can be considered. 
If the cost function $c$ involves the agent's state $s$ and action $a$, as well as the population distribution $L$, we say that the constraints are at the \textit{agent-population level}. Because now each agent's strategy does impact the value of the cost function $c$, when the mean-field flow $L$ is fixed, some strategies may violate the constraints and therefore the set of admissible strategies is only a subset of that of unconstrained MFGs. As a result, the Nash equilibria of CMFGs (as defined in Definition~\ref{def: constrained ne}) could be very different from the that of unconstrained MFGs. We shall discuss agent-population-level constraints in detail in Section~\ref{sec: agent-population-level}.

% When the cost functions have different forms, there are two types of constraints one can consider. When the cost function $c$ is independent of the agent's state $s$ and action $a$ and is only a function of the mean-field term $L$, the constraints are said to be at the \textit{population level}. As a consequence, each single agent's strategy does not impact whether the constraints are satisfied or not. This further means that CMFG NEs defined by Definition~\ref{def: constrained ne} are equivalent to MFG NEs with certain patterns. 

%{\color{blue}XXX.high-level examples, when this problem is important}

On the other hand, for a special case where the cost function $c$ depends only on the mean-field term $L$ and is independent of an agent's state $s$ and action $a$, the constraints are said to be at the \textit{population level}. Consequently, the strategy of an individual agent does not influence whether the constraints are satisfied. This further implies that CMFG Nash equilibria  are equivalent to MFG Nash equilibria with certain patterns. We refer to Section~\ref{sec: pop-level} for a more detailed discussion of population-level constraints.

To better illustrate the above ideas and to show the difference between population-level and agent-population-level constraints, we now introduce a concrete example of CMFGs, which is adapted from the Susceptible-Infected-Susceptible (SIS) model in \cite[Appendix A.3]{cui}. 

\begin{example}[Constrained SIS model] \label{exam: SIS}
The constrained SIS model evolves over discrete time steps ${0, 1, \dots, T}$ with finite horizon $T < \infty$. Each agent has two possible states, ${S, I}$, denoting ``susceptible'' and ``infected'', and two actions, ${U, D}$, representing ``going out'' and ``social distancing''. Agents incur a loss when infected but no penalty when susceptible. Social distancing lowers the probability of infection but yields a smaller immediate reward.

The reward function is defined as:
\begin{equation*}
    r_t(s_t, a_t, L_t) := -\delta_{s_t = I} - 0.5 \delta_{a_t = D},
\end{equation*}
where $\delta_\cdot$ is the indicator function, taking a value of 1 if the condition is true and 0 otherwise. The transition probabilities are:
\begin{align*}
    p_t(s_{t+1} = I \,|\, s_t = S, a_t = D, L_t) &= 0, \\
    p_t(s_{t+1} = I \,|\, s_t = S, a_t = U, L_t) &= 0.9 L_t(s_t = I, a_t = U), \\
    p_t(s_{t+1} = S \,|\, s_t = I, a_t = D, L_t) &= 0.5, \\
    p_t(s_{t+1} = S \,|\, s_t = I, a_t = U, L_t) &= 0.5 \Big[1 - 0.9 L_t(s_t = I, a_t = U)\Big].
\end{align*}
Here, $L_t(s_t = I, a_t = U) \in [0, 1]$ denotes the proportion of infected agents going out at time $t$.

To introduce constraints, we propose the following three illustrative examples: % (see Remark~\ref{rmk: macro_micro}):

\begin{align}
    \text{Agent-population-level (state):} \quad & \frac{1}{|\mathcal{T}|} \mathbb{E} \left( \sum_{t \in \mathcal{T}} \delta_{s_t = I} \,\bigg|\, s_0 \sim \mu_0 \right) \leq \gamma_0, \label{eq: micro_state} \\
    \text{Agent-population-level (action):} \quad & \frac{1}{|\mathcal{T}|} \mathbb{E} \left( \sum_{t \in \mathcal{T}} \delta_{a_t = U} \,\bigg|\, s_0 \sim \mu_0 \right) \leq \gamma_0, \label{eq: micro_action} \\
    \text{Population-level (state):} \quad & \frac{1}{|\mathcal{T}|} \sum_{t \in \mathcal{T}} L_t(s_t = I) \leq \gamma_0. \label{eq: macro_state}
\end{align}
Constraint \eqref{eq: micro_state} operates at the agent-population level, requiring that an individual agent's strategy ensures their infection risk remains below the threshold $\gamma_0$. This reflects a personalized safety measure that prioritizes the health of each agent. Similarly, constraint \eqref{eq: micro_action} is an agent-population-level restriction, limiting the time an individual agent spends outside. This could represent a government-imposed measure to control exposure during a pandemic. In contrast, constraint \eqref{eq: macro_state} is a population-level restriction. It does not enforce specific behavior for any individual agent. Instead, it sets a global limit on the infection rate in the population as a whole. This type of constraint aims to achieve collective safety and health outcomes by controlling the aggregate behavior of the agents.

%It does not enforce specific behavior for any individual agent. Instead, it sets a global limit on the infection rate in the population as a whole. Policies of all agents collectively need to satisfy this constraint.

%These constraints capture different social objectives. For example, constraint \eqref{eq: micro_state} is imposed on the agent-population-level, which says that the individual agent's policy has to guarantee that the agent's infection risk is lower than the threshold $\gamma_0$.  Constraint \eqref{eq: micro_action} is another constraint imposed on the agent-population-level,  limiting the individual agent's time of going out, which could be a reasonable constraint imposed by the government during pandemic.  On the other hand, constraint \eqref{eq: macro_state} is a population-level constraint. It does not restrict the policy of any single agent, but manages the overall infection rate across the population. 

The above description provides an intuitive explanation of the two types of constraints. In the following, we demonstrate through a concrete calculation that population-level constraints and agent-population-level constraints are indeed mathematically distinct.

For simplicity, let us consider the case where $T = 1$ and the initial distribution is $\mu_0(s_0 = I) = 1$. Denote by $L := L_0(s_0 = I, a_0 = U) \in [0, 1]$ the percentage of infected agents who choose to go out at  time 0, and by $\alpha := \mathbb{P}(a_0 = U) \in [0, 1]$ the probability that a specific agent chooses to go out given the population distribution $L$. A straightforward calculation shows that the expected reward for this specific agent is (we do not consider the action cost at $t=1$):
\[
\E (r_0 + r_1) = 0.5(1 - 0.9L) \alpha - 2.
\]
Since $1 - 0.9L > 0$, the agent maximizes their expected cumulative reward by always going out, \textit{i.e.}, $\alpha^\ast = 1$. This decision is independent of the population distribution $L$. Consequently, in the NE of the unconstrained MFG, the percentage of infected agents at $t = 1$ will be:
\[
\mathbb{P}(s_1 = I) = 0.95.
\]
As a comparison, if all agents instead choose social distancing at the initial time, only half of them will remain infected at $t = 1$.

Now, suppose we impose the following \textbf{population-level constraint}:
\begin{equation} \label{eq: example_pop}
    L_1(s_1 = I) \leq \gamma_0.
\end{equation}
Based on the earlier calculation, there will be no NE in the CMFG if $\gamma_0 < 0.95$. Alternatively, suppose we instead impose the following \textbf{agent-population-level constraint}:
\begin{equation} \label{eq: example_agent}
    \E(\delta_{s_1 = I}) \leq \gamma_0.
\end{equation}
Under the constraint \eqref{eq: example_agent}, there are three cases to discuss:
\begin{itemize}
    \item Case 1: $\gamma_0 > 0.5$

    In this case, there will be an NE in the CMFG. To spell out, by Definition~\ref{def: constrained ne}, the NE policy $\alpha^{\rm NE}$ solves the following system of equations:
    \begin{align*}
        \begin{cases}
            0.5 + 0.45 \alpha^{\rm NE} L^{\rm NE} = \min\{\gamma_0,\, 0.95\}, \\
            \alpha^{\rm NE} = L^{\rm NE},
        \end{cases}
    \end{align*}
    which gives:
    \begin{equation*}
        \alpha^{\rm NE} = \min \left\{\sqrt{\frac{\gamma_0 - 0.5}{0.45}},\, 1\right\}.
    \end{equation*}

    \item Case 2: $\gamma_0 = 0.5$

    In this case, there will be no NE in the CMFG. Since $\E(\delta_{s_1 = I}) = 0.5 + 0.45 L \alpha$, the constraint translates into $\alpha L = 0$. When $L = 0$, every $\alpha \in [0, 1]$ is feasible. However, the optimal action $\alpha^\ast = 1$, which contradicts the consistency condition of the NE. On the other hand, when $L \in (0, 1]$, the only feasible action is $\alpha = 0$, which contradicts the consistency condition, too. As a consequence, there will be no NE.

    \item Case 3: $\gamma_0 < 0.5$

    In this case, there will be no NE in the CMFG. This is obvious since there will be no feasible action.
\end{itemize}
This concrete calculation thus demonstrates the difference between population-level constraints and agent-population-level constraints.
$\hfill\square$
\end{example}

\begin{remark}[Lagrangian method and penalty method]
Alternative approaches to handling constraints include the Lagrangian method and penalty method. For Lagrangian method, the representative agent solves the following unconstrained MDP problem:
\begin{align}
    \text{maximize}_{\{ \pi_t \}_{t \in \T}} & \quad \E \left[ \sum_{t \in \T} \Big(r_t(s_t, a_t, L_t) - \eta_1^\top c_t(s_t, a_t, L_t)\Big) \,\bigg|\, s_0 \sim \mu_0 \right] \nonumber \\
    \text{subject to} & \quad a_t \sim \pi_t(s_t), \quad s_{t + 1} \sim p_t(s_t, a_t, L_t), \nonumber
\end{align}
where $\eta_1 > 0$ is the Lagrangian multiplier. 

For penalty method, the representative agent solves the following unconstrained MDP problem:
\begin{align}
    \text{maximize}_{\{ \pi_t \}_{t \in \T}} & \quad \E \left[ \sum_{t \in \T} \Big(r_t(s_t, a_t, L_t) - \eta_2^\top (c_t(s_t, a_t, L_t)-\gamma_0)^+\Big) \,\bigg|\, s_0 \sim \mu_0 \right] \nonumber \\
    \text{subject to} & \quad a_t \sim \pi_t(s_t), \quad s_{t + 1} \sim p_t(s_t, a_t, L_t), \nonumber
\end{align}
where $\eta_2 > 0$ is the penalty coefficient. For these two problems, the unconstrained MDPs always admit  solutions, and the corresponding MFGs fall into the standard setting. 
However, studying the constrained version directly is often preferable or necessary, especially when violating the constraints may lead to detrimental consequences and hard constraints must be enforced. For Lagrangian method, the relationship between $\gamma_0$ in \eqref{eq: cmdp3} and the Lagrangian multiplier coefficient $\eta_1$ is unclear, making it challenging to define an appropriate $\eta_1$ in practice. And for penalty method, one needs to increase the value of $\eta_2$ to make the problem asymptotically equivalent to the original constrained problem.

%as constraints appear in the form of \eqref{eq: cmdp3} in most real-world scenarios, and 
\end{remark}

\section{Mean-Field Games with Agent-Population-Level Constraints} \label{sec: agent-population-level}

In this section, we focus on the problem of CMFGs with general agent-population-level constraints, \textit{i.e.}, the cost function $c$ involves agent's state $s$ and action $a$, as well as population distribution $L$. We first analyze the existence of CMFG NEs defined by Definition~\ref{def: constrained ne}. We then provide a numerical scheme to solve the CMFGs.

\subsection{Existence and Uniqueness of Nash Equilibrium in CMFGs} \label{sec: properties}

Our main result in this section is that a CMFG NE defined by Definition~\ref{def: constrained ne} always exists under proper continuity assumptions and a feasibility assumption on the constraints of the CMFG.

The first assumption in the following ensures that for any individual constraint indexed by $i$, there exists a policy $\pi^{(i)}$ that satisfies all constraints while being strictly feasible for constraint $i$.  

\begin{assumption}[Strict feasibility] \label{ass: strict feasibility}
    There exists $\delta > 0$ such that for any mean-field flow $L$ and any index $1 \leq i \leq k$, there exists a policy $\pi^{(i)}$ satisfying:
    \begin{equation*}
        \mathbb{E}_{a \sim \pi^{(i)}} \left[ \sum_{t \in \mathcal{T}} c_t (s_t, a_t, L_t) \,\Big|\, s_0 \sim \mu_0 \right] \leq \gamma_0,
    \end{equation*}
    and
    \begin{equation*}
        \mathbb{E}_{a \sim \pi^{(i)}} \left[ \sum_{t \in \mathcal{T}} [c_t]_i (s_t, a_t, L_t) \,\Big|\, s_0 \sim \mu_0 \right] \leq [\gamma_0]_i - \delta,
    \end{equation*}
    where $[\cdot]_i$ denotes the $i$-th entry of the vector.
\end{assumption}

When $\delta=0$, the assumption above simply ensures the feasibility of the constraints, guaranteeing that for any mean-field flow, at least one policy satisfies the given conditions. This can also be interpreted as stating that each agent always has at least one feasible strategy, regardless of the strategies chosen by competitors. By assuming $\delta > 0$, we strengthen this feasibility condition by ensuring the existence of a \textit{strictly feasible} policy for each constraint. As we will show, this tightening can not be relaxed and plays a crucial role in our subsequent analysis.

% {\color{red}XXX. Add a comment that the intuition is that we have agent-population-level feasibility facing arbitrary/adversarial competitors (thus for all $L$ stuff in the assumption). This explanation may make the assumption more naturally understandable from the original N-player perspectives.}

% We incorporate mean-field flow $L$ into agent-population-level constraints since the agents have mean-field-type interactions (\textit{e.g.}, competition). In the limiting case where $\delta = 0$, Assumption~\ref{ass: strict feasibility} 

We also impose the following continuity assumption.
Such continuity conditions on rewards and transition kernels are standard in the existing  MFG literature to guarantee the existence of NEs. Here, we extend this assumption to include the continuity of the constraints as well.

\begin{assumption}[Continuity] \label{ass: continuous}
    The functions $p$, $r$, and $c$ are  continuous with respect to the mean-field flow. %That is, there exist constants $C_p, C_r, C_c > 0$, such that for any $t \in \mathcal{T}$ and any $L^1,L^2\in \Delta(\mathcal{S}\times\mathcal{A})$, we have the inequalities:
    % \begin{align*}
    %     \left|\left|\sum_{s \in \mathcal{S}, a \in \mathcal{A}} \Big|p_t\left(s, a, L^1\right) - p_t\left(s, a, L^2\right)\Big| \right|\right|_1 & \leq C_p \left|\left|L^1 - L^2\right|\right|_1, \\ 
    %     \left|\left|r_t\left(L^1\right) - r_t\left(L^2\right)\right|\right|_1 & \leq C_r \left|\left|L^1 - L^2\right|\right|_1, \\
    %     \left|\left|\sum_{s \in \mathcal{S}, a \in \mathcal{A}} \Big|c_t\left(s, a, L^1\right) - c_t\left(s, a, L^2\right)\Big| \right|\right|_1 & \leq C_c \left|\left|L^1 - L^2\right|\right|_1,
    % \end{align*}
    % where $p_t\left(s, a, L_t\right)$ is treated as a vector in $\mathbb{R}^{|\mathcal{S}|}$, $r_t(L_t)$ and $L_t$ are treated as vectors in $\mathbb{R}^{|\mathcal{S}| |\mathcal{A}|}$, and $c_t(s, a, L_t)$ is treated as a vector in $\mathbb{R}^{k}$. 
\end{assumption}

Under these assumptions, we are able to establish the existence of NE for the constrained MFGs.
\begin{theorem}[Existence of NE in CMFGs] \label{thm: exist}
    Under Assumptions~\ref{ass: strict feasibility} and \ref{ass: continuous}, there exists at least one CMFG NE defined by Definition~\ref{def: constrained ne}.
\end{theorem}

\begin{remark}
    Assumption~\ref{ass: strict feasibility} is critical for the proof of existence. Indeed, the CMFG may not have an NE when Assumption~\ref{ass: strict feasibility} is not satisfied. For instance, let us set $\gamma_0 = 0.5$ in \eqref{eq: example_agent} (\textit{cf.} Example~\ref{exam: SIS}). Then there is no positive $\delta$ satisfying Assumption~\ref{ass: strict feasibility}. By the analysis in Example~\ref{exam: SIS}, the resulting CMFG  %$\alpha=0$ is a feasible agent policy for any given population distribution $L \in [0, 1]$. However, 
    does not have any NE. We further notice that the CMFG also satisfies Assumption~\ref{ass: continuous}, which means Assumption~\ref{ass: strict feasibility} cannot be relaxed.
\end{remark}

The proof of Theorem \ref{thm: exist} relies on the following two lemmas. Lemma~\ref{lemma: lp} reformulates the constrained Markov Decision Process (MDP) problem as a linear program (LP). Lemma~\ref{lemma: dual selection} then establishes that the solution to the Karush-Kuhn-Tucker (KKT) conditions, arising from the LP formulation in Lemma~\ref{lemma: lp}, can always be found within a bounded domain, where the bound is independent of the choice of the mean-field flow.

To aid our later discussion, we first introduce the concept of \textit{occupation measure}, which characterizes the distribution of state-action of a single player. Specifically, when the mean-field flow $L$ is fixed, for any admissible policy $\pi$, its associated occupation measure $d \in \big[\Delta(\state \times \action)\big]^{|\T|}$ describes the probability of the representative player staying at state $s$ and taking action $a$ at time $t$. %Specifically, one can derive the following recursive equations for calculating the occupation measure.  
With a slight abuse of notation, for any mean-field flow $L$, $s \in \state$ and $a \in \action$, define $\Psi(\pi, L) \in \big[\Delta(\state \times \action)\big]^{|\T|}$ recursively by:
\begin{align*}
    \begin{cases}
         \Psi_0(\pi, L) (s, a) := \pi_0(a \,|\, s) \mu_0(s), \\
        \Psi_{t + 1}(\pi, L) (s, a) := \pi_{t+1}(a\,|\,s)\sum_{s'\in\mathcal{S},a'\in\mathcal{A}}p_t(s \,|\, s',a',L_t)\Psi_{t}(\pi, L)(s',a'). 
    \end{cases}
\end{align*}
Then $\Psi(\pi,L)$ is the occupation measure of a representative agent taking policy $\pi$ under mean-field flow $L$.

%\footnote{Occupation measure describes the distribution of a specific agent. Note its conceptual difference with a mean-field flow (the latter characterizes the distribution of a population).} 
On the other hand, for any $d \in \big[\Delta(\state \times \action)\big]^{|\T|}$, we can retrieve the policy which leads to the occupation measure as follows. Specifically, define $\Pi(d)$ by the set of policies such that:
\begin{equation*}
    \forall t \in \T, s \in \state, a \in \action, \quad \pi_t(a \,|\, s) = 
    \begin{cases}
        d_t(s, a) \,\Big\slash\, \sum_{a^\prime \in \action} d_t(s, a^\prime), \quad \text{if $\sum_{a^\prime \in \action} d_t(s, a^\prime)$ > 0}, \\
        \text{any vector in $\Delta(\action)$}, \quad \text{else}.
    \end{cases}
\end{equation*}
Then any policy $\pi\in\Pi(d)$  induces the occupation measure $d$.

Utilizing the connection between admissible policies and occupation measures, one can transform the constrained MDP problem into a linear program on occupation measure.
\begin{lemma}[CMDPs as LPs] \label{lemma: lp}    Given the mean-field flow $L$, solving the constrained MDP problem $\cmdp{L}$ is equivalent to solving the following LP problem:
    \begin{align}
        \text{minimize}_{d} & \quad -r_L^\top d \label{eq: lp1} \tag{LP 1} \\
        \text{subject to} & \quad A_L d = b, \quad c_L d \leq \gamma_0, \quad d \geq 0. \label{eq: lp2} \tag{LP 2}
    \end{align}
    Here $r_L \in \R^{|\state||\action||\T|}$, $b \in \R^{|\state||\T|}$, $c_L \in \R^{k \times |\state||\action||\T|}$, and $A_L \in \R^{|\state||\T| \times |\state||\action||\T|}$ are defined by:
    \begin{align*}
        r_L & := 
        \begin{pmatrix}
            r_0 (L_0) \\
            r_1 (L_1) \\
            \vdots \\
            r_T (L_T) \\
        \end{pmatrix}, 
        \quad r_t (L_t) \in \R^{|\state||\action|}, \quad b := 
        \begin{pmatrix}
            0 \\
            0 \\
            \vdots \\
            0 \\
            \mu_0
        \end{pmatrix}, 
        \quad \mu_0 \in \R^{|\state|}, \\
        c_L & := 
        \begin{pmatrix}
            \Big([c_0]_1 (L_0)\Big)^\top & \Big([c_1]_1 (L_1)\Big)^\top & \cdots & \Big([c_T]_1 (L_T)\Big)^\top \\
            \Big([c_0]_2 (L_0)\Big)^\top & \Big([c_1]_2 (L_1)\Big)^\top & \cdots & \Big([c_T]_2 (L_T)\Big)^\top \\
            \vdots & \vdots & \ddots & \vdots \\
            \Big([c_0]_k (L_0)\Big)^\top & \Big([c_1]_k (L_1)\Big)^\top & \cdots & \Big([c_T]_k (L_T)\Big)^\top
        \end{pmatrix}, \quad [c_t]_i (L_t) \in \R^{|\state||\action|}, \\
        A_L & := 
        \begin{pmatrix}
            A_0 (L_0) & -Z & 0 & 0 & \cdots & 0 & 0 \\
            0 & A_1 (L_1) & -Z & 0 & \cdots & 0 & 0 \\
            \vdots & \vdots & \vdots & \vdots & \ddots & \vdots & \vdots \\
            0 & 0 & 0 & 0 & \cdots & A_{T-1} (L_{T-1}) & -Z \\
            Z & 0 & 0 & 0 & \cdots & 0 & 0
        \end{pmatrix}, \quad Z := \overbrace{\Big(I_{|\state|}, I_{|\state|}, \cdots, I_{|\state|} \Big)}^{|\action|}, 
    \end{align*}
    where $I_{|\state|}$ denotes the identity matrix of shape $|\state| \times |\state|$, and for each $t \in \T\backslash\{T\}$, $A_t(L_t) \in \R^{|\state| \times |\state||\action|}$ is defined by:
    \begin{align*}
        \forall 1 \leq i \leq |\state|, \quad \text{the $i$-th row of } A_t(L_t) := \Big(p_t (i \,|\, \cdot, \cdot, L_t)\Big)^\top \in \R^{1 \times |\state||\action|}.
    \end{align*}
    Specifically, if $\pi^\ast$ is an optimal policy of $\cmdp{L}$, then there exists $d^\ast$ such that $\pi^\ast \in \Pi(d^\ast)$ and $d^\ast$ is an optimal solution of the above LP problem. Conversely, if $d^\ast$ is an optimal solution of the above LP problem, then for any $\pi^\ast \in \Pi(d^\ast)$, $\pi^\ast$ is an optimal policy of $\cmdp{L}$.
\end{lemma}

The LP formulation of MDPs has been exploited in \cite{mfomo}, and the proof of Lemma~\ref{lemma: lp} follows from that of \cite[Lemma 2, Lemma 3]{mfomo} and thus is omitted here.

By Lemma~\ref{lemma: lp} and the strong duality of LP (see \cite{luenberger1984linear}, for example), solving $\cmdp{L}$ is then equivalent to finding the solution that satisfies the following KKT conditions.

\begin{kkt} \label{kkt}
    \begin{align*}
        \text{find} \quad & \left\{d \in \R^{|\state||\action||\T|}, y \in \R^{|\state||\T|}, z \in \R^{|\state||\action||\T|}, \lambda \in \R^k \right\} \nonumber \\
        \text{such that} \quad & -r_L = A_L^\top y + z - c_L^\top \lambda, \quad A_L d = b, \quad c_L d \leq \gamma_0, \quad d \geq 0, \\
        & z \geq 0, \quad z^\top d = 0, \quad \lambda \geq 0, \quad \lambda^\top (c_L d - \gamma_0) = 0. 
    \end{align*}
\end{kkt}
\noindent
Since the KKT conditions are uniquely determined by the mean-field flow $L$, we denote them by $\K(L)$ in the rest of our paper. It is obvious that finding an NE in a CMFG is equivalent to finding a fixed point of the KKT conditions, in the following sense.

\begin{proposition} \label{prop: fixed_point}
    $(\pi, L)$ is a CMFG NE defined by Definition~\ref{def: constrained ne} if and only if $\pi \in \Pi(L)$ and there exists $(y, z, \lambda)$ such that $(L, y, z, \lambda)$ is a solution of $\K(L)$.
\end{proposition}

The following lemma guarantees that one can always find solutions $y,z,\lambda$ in a compact set.
\begin{lemma}[A bounded solution of $\K(L)$] \label{lemma: dual selection}
    Given the mean-field flow $L$, under Assumption~\ref{ass: strict feasibility}, let $(d, y, z, \lambda)$ be a solution of $\K(L)$. Then, $\left|\left|\lambda\right|\right|_\infty \leq |\T| \rmax \slash \delta$, and there exists $\left(y^\prime, z^\prime\right)$ in a bounded domain such that $\left(d, y^\prime, z^\prime, \lambda\right)$ is a solution of $\K(L)$. Specifically, the bounds are given by:
    \begin{align*}
        \left|\left|y^\prime\right|\right|_1 & \leq \frac{|\mathcal{S}| |\T| (|\T| + 1)}{2} \rmax \left(1 + \frac{|\T|}{\delta} \cmax \right), \\
        \left|\left|z^\prime\right|\right|_1 & \leq |\mathcal{S}| |\mathcal{A}| (|\T|^2 - |\T| + 2) \rmax \left(1 + \frac{|\T|}{\delta} \cmax \right),
    \end{align*}
    where $\rmax$ and $\cmax$ denote the maximum values of $r$ and $c$, respectively, \textnormal{i.e.}, 
    \begin{align*}
        \rmax & := \max_{t \in \T, s \in \state, a \in \action, L \in \Delta(\state \times \action)} r_t(s, a, L), \\
        \cmax & := \max_{t \in \T, s \in \state, a \in \action, L \in \Delta(\state \times \action)} ||c_t(s, a, L)||_\infty.
    \end{align*}
\end{lemma}

\begin{proof}
    Suppose $(d, y, z, \lambda)$ is a solution of $\K(L)$. We first prove that $||\lambda||_\infty \leq |\T| r_{\max} \slash \delta$. To see this, notice that $d$ is a minimizer of the LP problem:
    \begin{align*}
        \text{minimize}_{d} & \quad -\left(r_L - c_L^\top \lambda \right)^\top d \\
        \text{subject to} & \quad A_L d = b, \quad d \geq 0.
    \end{align*}
    Therefore, for any index $1 \leq i \leq k$, choosing policy $\pi^{(i)}$ defined in Assumption~\ref{ass: strict feasibility} and defining the associated occupation measure $d^{(i)} := \Psi\left(\pi^{(i)}, L\right)$, we will have:
    \begin{align}
        \left(r_L - c_L^\top \lambda\right)^\top d & \geq \left(r_L - c_L^\top \lambda\right)^\top d^{(i)}, \nonumber \\
        r_L^\top \left(d - d^{(i)}\right) & \geq \lambda^\top c_L \left(d - d^{(i)}\right). \label{eq: bounded kkt}
    \end{align}
    In \eqref{eq: bounded kkt}, LHS $\leq |\T|\rmax$, and RHS $= \lambda^\top \left(\gamma_0 - c_L d^{(i)}\right) \geq [\lambda]_i \left[\gamma_0 - c_L d^{(i)}\right]_i \geq \delta [\lambda]_i$. Thus, we see that $[\lambda]_i \leq |\T| r_{\max} \slash \delta$, which yields the desired bound of $\lambda$.

    The boundedness of $y^\prime$ and $z^\prime$ is directly implied by combining the proof of \cite[Proposition 6]{mfomo} and . The only difference with \cite[Proposition 6]{mfomo} is that $r_L$ is replaced by $r_L - c_L^\top \lambda$. This finalizes our proof.
\end{proof}

Finally, we present the proof of Theorem~\ref{thm: exist} below.

\begin{proof}[Proof of Theorem~\ref{thm: exist}]
    We start with defining the set-valued mapping $\Gamma$, which maps a mean-field flow $L$ to the set of optimal solutions  of the LP introduced in Lemma~\ref{lemma: lp}. By Proposition~\ref{prop: fixed_point}, it suffices to establish that $\Gamma$ admits a fixed point. In the following, we prove this by aplying Kakutani's fixed-point theorem. 
    
    To proceed, we verify  the conditions of Kakutani's fixed-point theorem. First, notice that the set of mean-field flows can be viewed as a concatenation of $|\T|$ probability simplices, each of dimension $|\state||\action| - 1$, making it a non-empty, compact and convex subset of $\R^{|\T| (|\state||\action| - 1)}$.
    Second, the convexity of $\Gamma(L)$  for any mean-field flow $L$ follows from the linearity  of \eqref{eq: lp1} and \eqref{eq: lp2}. 
    Finally, we verify that $\Gamma$ has a closed graph. Suppose $\{L_n\}_{n \geq 0}$ and $\{d_n\}_{n \geq 0}$ are  sequences of mean-field flows and occupation measures such that 
    \begin{enumerate}
        \item $d_n \in \Gamma(L_n)$ for all $n$, and
        \item  there exist $\hat L$ and $\hat d$ such that $\lim_{n \rightarrow \infty} L_n = \hat{L}$ and $\lim_{n \rightarrow \infty} d_n = \hat{d}$.
    \end{enumerate}
    To show that $\hat{d} \in \Gamma(\hat{L})$, we note that for each $n$, there exists $(y_n, z_n, \lambda_n)$ such that $(d_n, y_n, z_n, \lambda_n)$ is a solution of $\K(L_n)$. By Lemma~\ref{lemma: dual selection}, the sequence $\{ (d_n, y_n, z_n, \lambda_n) \}_{n \geq 0}$ is uniformly bounded, ensuring the existence of a convergent subsequence. Let $(\hat{d}, \hat{y}, \hat{z}, \hat{\lambda})$ be its limit. Then by Assumption~\ref{ass: continuous}, this limit is a solution to $\K(\hat{L})$, which implies $\hat{d} \in \Gamma(\hat{L})$. This completes the proof.
    % \emph{1)} $d_n \in \Gamma(L_n)$ for all $n$ and that \emph{2)} With the aim of 
\end{proof}

\begin{remark}
    The proof of Theorem~\ref{thm: exist} differs from  existing methods to prove the existence of NEs, and highly relies on the primal-dual formulation. Specifically, the fixed point mapping in our approach is defined in terms of mean-field flows (and occupation measures), with the analysis depending on the boundedness of the dual variables (Lemma~\ref{lemma: dual selection}). In contrast, existing approaches, such as those in \cite{cui}, define the fixed point mapping via optimal policies and analyze it using value functions and Q-functions. However, this approach is difficult to apply to CMFGs due to the absence of Bellman equations.  % in  We mention that their method fails for CMFGs as the value functions are not well-defined in CMFGs. The core of our proof is Lemma~\ref{lemma: dual selection}, \textnormal{i.e.}, the dual variables can be selected in a bounded domain.
\end{remark}

% {\color{red}
% \begin{remark}
%     XXX. Here the proof to show the existence of NE to the CMFG  is different from most of the existing methods to prove the existence of NE for MFGs, where the fixed point mapping is defined via policies and mean-field flows. Comment on the difference and highlight that our method fully utilize the occupation measures. ($\lambda$ bounded)
% \end{remark}
% }

Beyond existence, we also establish a uniqueness result under the classical Lasry-Lions monotonicity condition.

\begin{assumption}[Monotonicity condition]\label{ass:monotone}
    The transition probabilities and the constraints are independent of the population distribution. In addition, the reward functions satisfy the following monotonicity condition: for any two mean-field flows $L^1$ and $L^2$,
    \begin{equation*}
        \sum_{t\in\mathcal{T}}\sum_{s\in\mathcal{S}}\sum_{a\in\mathcal{A}} \big(r_t(s,a,L^1_t)-r_t(s,a,L^2_t)\big)\big(L_t^1(s,a)-L_t^2(s,a)\big)\leq 0,
    \end{equation*}
    with equality if and only if $L^1 \equiv L^2$.
\end{assumption}

\begin{theorem}[Uniqueness of CMFG NE]
    Under Assumptions \ref{ass: strict feasibility}, \ref{ass: continuous}, and \ref{ass:monotone}, the CMFG admits a unique Nash equilibrium mean-field flow as defined in Definition \ref{def: constrained ne}.
\end{theorem}

\begin{proof}
    By Theorem \ref{thm: exist}, there exists at least one CMFG NE. Suppose for contradiction that there are two distinct NE mean-field flows, $L^1$ and $L^2$. Then, by the definition of CMFG NE and Lemma \ref{lemma: lp},
    \begin{align*}
        r_{L^1}^\top L^1 \geq r_{L^1}^\top L^2, \quad 
        r_{L^2}^\top L^2 \geq r_{L^2}^\top L^1.
    \end{align*}
    This yields
    \begin{align*}
        0> (r_{L^1}-r_{L^2})^\top (L^1-L^2)=r_{L^1}^\top L^1 - r_{L^1}^\top L^2 + r_{L^2}^\top L^2 - r_{L^2}^\top L^1 \geq 0,
    \end{align*}
a contradiction. Hence the equilibrium is unique in terms of the mean-field flow.
\end{proof}

\subsection{Optimization Framework for Solving NE in CMFGs} \label{sec: cmfomo}

In this section, we provide a framework of \textit{Constrained Mean-Field Occupation Measure Optimization} (CMFOMO), which transforms the problem of solving NEs 
%fixed point problem (\textit{cf.} Proposition~\ref{prop: fixed_point}) 
into an optimization problem with bounded variables and convex constraints. This framework extends \cite[Theorem 5]{mfomo} to incorporate constraints in the representative agent's problem. 
The following Theorem~\ref{thm: cmfomo} is a direct implication of Definition~\ref{def: constrained ne} and Lemmas~\ref{lemma: lp} and \ref{lemma: dual selection}.

\begin{theorem}[Constrained MF-OMO] \label{thm: cmfomo}
    Let $c_i > 0$ $(i = 1, 2, 3, 4, 5)$. Under Assumption~\ref{ass: strict feasibility}, $(\pi^\ast, L^\ast)$ is a CMFG NE defined by Definition~\ref{def: constrained ne} if and only if there exists $(y^\ast, z^\ast, \lambda^\ast, w^\ast)$ such that $(L^\ast, y^\ast, z^\ast, \lambda^\ast, w^\ast)$ is an optimal solution of the following optimization problem with the optimal objective value being 0: %{\color{red} isn't it supposed to be 0? Assumption 1 guarantees the existence}
    \begin{align}
        \text{minimize}_{L, y, z, \lambda, w} \quad & c_1 \left|\left| A_L^\top y + z + r_L - c_L^\top \lambda\right|\right|_2 + c_2 \left|\left|A_L L - b\right|\right|_2 + c_3 (z^\top L) \nonumber \\
        & \qquad \qquad \qquad \qquad \qquad \quad \;\, + c_4 \left|\left|\gamma_0 - c_L L - w\right|\right|_2 + c_5 \left|\lambda^\top (\gamma_0 - c_L L)\right| \tag{CMFOMO} \label{eq: cmfomo obj} \\
        \text{subject to} \quad & L \geq 0, \quad \rv{1}^\top L_t = 1, \quad \forall t \in \T, \nonumber \\
        & ||y||_1 \leq \frac{|\mathcal{S}| |\T| (|\T| + 1)}{2} r_{\max} \left(1 + \frac{|\T|}{\delta} c_{\max}\right), \nonumber \\
        & z \geq 0, \quad ||z||_1 \leq |\mathcal{S}| |\mathcal{A}| \left(|\T|^2 - |\T| + 2\right) r_{\max} \left(1 + \frac{|\T|}{\delta} c_{\max}\right), \nonumber \\
        & 0 \leq \lambda \leq \frac{|\T| r_{\max}}{\delta}, \nonumber \\
        & 0 \leq w \leq \gamma_0, \nonumber
    \end{align}
    where the mean-field flow $L$ is treated as a vector in $\R^{|\state||\action||\T|}$, $L_t$ is defined by the slice of $L$ whose index ranges from $(|\state||\action|t + 1)$ to $|\state||\action|(t+1)$, and $\rv{1}$ is defined by a vector of length $|\state||\action|$ whose entries are equal to 1. We refer to Assumption~\ref{ass: strict feasibility} and Lemma~\ref{lemma: lp} for the definition of the rest of the notations.
\end{theorem}

\begin{remark}
    In the objective of CMFOMO, $c_2$ and $c_4$ control the consistency (\textit{cf.} Remark~\ref{rmk: opt_consis}) and feasibility of $L$, respectively; $c_1$ characterizes the optimality condition (\textit{cf.} Remark~\ref{rmk: opt_consis}); $c_3$ and $c_5$ are the complementarity terms connecting the previous conditions into a single optimization problem. When $c_4 = c_5 = 0$ and fixing $\lambda = 0$, the objective of CMFOMO reduces to MF-OMO introduced in \cite{mfomo}. 
\end{remark}

When allowing optimization errors, as we shall show next, $\epsilon$-optimal solutions of CMFOMO are approximations of the NEs of the CMFG. To quantify the gaps between approximated NEs and true NEs, we introduce the following metrics.

\begin{definition}[Optimality gap] \label{def: opt_gap}
    For any policy $\pi$, we define its \textit{optimality gap} by:
    \begin{equation*}
        G_{\text{opt}}(\pi) := V^\ast_c (\Psi(\pi)) - V^\pi (\Psi(\pi)),
    \end{equation*}
    where $V^\ast_c (\Psi(\pi))$ denotes the maximum expected cumulative reward of $\cmdp{\Psi(\pi)}$ (which always exists under Assumption~\ref{ass: strict feasibility}), and $V^\pi (\Psi(\pi))$ denotes the expected cumulative reward of policy $\pi$ under the mean-field flow $\Psi(\pi)$, \textnormal{i.e.},
    \begin{equation*}
        V^\pi(\Psi(\pi)) := \mathbb{E}_{a \sim \pi} \left[ \sum_{t \in \mathcal{T}} r_t(s_t, a_t, \Psi_t(\pi)) \,\bigg|\, s_0 \sim \mu_0 \right].
    \end{equation*}
\end{definition}

\begin{definition}[Feasibility gap] \label{def: fea_gap}
    For any policy $\pi$ and mean-field flow $L$, we define the \textit{feasibility gap} of $\pi$ under $L$ by:
    \begin{equation*}
        G_{\text{fea}}(\pi, L) := \left|\left|\min\Big(0, \gamma_0 - \gamma_L(\pi) \Big)\right|\right|_2,
    \end{equation*}
    where the $\min$ operator is applied entry-wise, and
    \begin{equation*}
        \gamma_L (\pi) := \mathbb{E}_{a \sim \pi} \left[ \sum_{t \in \mathcal{T}} c_t(s_t, a_t, L_t) \,\bigg|\, s_0 \sim \mu_0 \right].
    \end{equation*}
    With an abuse of notation, we define the feasibility gap of $\pi$ by:
    \begin{equation*}
        G_{\text{fea}}(\pi) := G_{\text{fea}}(\pi, \Psi(\pi)).
    \end{equation*}
\end{definition}

\begin{remark}
    For any given policy $\pi$, $G_{\text{opt}}(\pi)$ measures its performance compared with the optimal policy when the population distribution is induced by $\pi$, while $G_{\text{fea}}(\pi)$ measures the extent to which $\pi$ violates the constraints. 
\end{remark}

When $G_{\text{fea}}(\pi, L) > 0$, we say that $\pi$ is \textit{infeasible} under $L$. Otherwise, we say $\pi$ is \textit{feasible} under $L$. If $\pi$ is feasible under $\Psi(\pi)$, then by definition, $G_{\text{opt}}(\pi)$ will always be non-negative. However, if $\pi$ is infeasible under $\Psi(\pi)$, $G_{\text{opt}}(\pi)$ may be negative.

In order to obtain the non-asymptotic analysis, we impose a stronger assumption in the subsequent analysis, which requires the rewards, transition kernels and constraints to be Lipschitz continuous in the mean-field flow. 
\begin{assumption}[Lipschitz continuity] \label{ass: Lip_continuous}
    The functions $p$, $r$, and $c$ are Lipschitz continuous with respect to the mean-field flow. That is, there exist constants $C_p, C_r, C_c > 0$, such that for any $t \in \mathcal{T}$ and any $L^1,L^2\in \Delta(\mathcal{S}\times\mathcal{A})$, we have the inequalities:
    \begin{align*}
        \left|\left|\sum_{s \in \mathcal{S}, a \in \mathcal{A}} \Big|p_t\left(s, a, L^1\right) - p_t\left(s, a, L^2\right)\Big| \right|\right|_1 & \leq C_p \left|\left|L^1 - L^2\right|\right|_1, \\ 
        \left|\left|r_t\left(L^1\right) - r_t\left(L^2\right)\right|\right|_1 & \leq C_r \left|\left|L^1 - L^2\right|\right|_1, \\
        \left|\left|\sum_{s \in \mathcal{S}, a \in \mathcal{A}} \Big|c_t\left(s, a, L^1\right) - c_t\left(s, a, L^2\right)\Big| \right|\right|_1 & \leq C_c \left|\left|L^1 - L^2\right|\right|_1,
    \end{align*}
    where $p_t\left(s, a, L_t\right)$ is treated as a vector in $\mathbb{R}^{|\mathcal{S}|}$, $r_t(L_t)$ and $L_t$ are treated as vectors in $\mathbb{R}^{|\mathcal{S}| |\mathcal{A}|}$, and $c_t(s, a, L_t)$ is treated as a vector in $\mathbb{R}^{k}$. 
\end{assumption}

The following theorem shows that any $\epsilon$-sub-optimal solution of CMFOMO leads to an approximated NE of the CMFG with explicit approximation error bounded by a linear function of $\epsilon$. Its proof can be found in Appendix~\ref{app: proofs}.

\begin{theorem}[Suboptimality of CMFOMO] \label{thm: suboptimal cmfomo}
    Let $\epsilon \geq 0$. Under Assumption \ref{ass: strict feasibility} and Assumption~\ref{ass: Lip_continuous}, suppose $(L, y, z, \lambda, w)$ is a feasible solution of CMFOMO defined in Theorem~\ref{thm: cmfomo} with its objective \eqref{eq: cmfomo obj} $\leq \epsilon$. Then, for any $\pi \in \Pi(L)$, its optimality gap and feasibility gap are bounded by:
    \begin{align}
        - \epsilon \sqrt{k} \frac{|\T| \rmax}{\delta} \zeta_3 \leq G_{\text{opt}}(\pi) & \leq \epsilon \Big[(|\T| + 1) \zeta_1 + \zeta_2 + \zeta_4 \Big], \label{eq: opt estimate} \\
        0 \leq G_{\text{fea}}(\pi) & \leq \epsilon \zeta_3, \label{eq: fea estimate}
    \end{align}
    where
    \begin{align}
        \beta_y & := \frac{|\mathcal{S}| |\T| (|\T| + 1)}{2} r_{\max} \left(1 + \frac{|\T|}{\delta} c_{\max}\right), \nonumber \\
        \beta_z & := |\mathcal{S}| |\mathcal{A}| (|\T|^2 - |\T| + 2) r_{\max} \left(1 + \frac{|\T|}{\delta} c_{\max}\right), \nonumber \\
        \alpha (c_2) & := \frac{1}{c_2} \frac{(C_p + 1)^{|\T| + 1} - 2C_p - 1}{C_p^2} \sqrt{|\mathcal{S}|}, \label{eq: alpha} \\
        \zeta_1(c_1, c_2) & := \frac{1}{c_1} + \alpha(c_2) \left( C_p \beta_y + C_r + k C_c \frac{|\T| r_{\max}}{\delta} \right), \nonumber \\
        \zeta_2(c_2, c_3) & := \frac{1}{c_3} + \beta_z \alpha(c_2), \nonumber \\
        \zeta_3(c_2, c_4) & := \frac{1}{c_4} + \alpha(c_2) \Big(k c_{\max}  + C_c |\mathcal{S}| |\mathcal{A}| |\T|\Big), \label{eq: zeta_3} \\
        \zeta_4(c_2, c_5) & := \frac{1}{c_5} + k \alpha(c_2) \frac{|\T| r_{\max}}{\delta} \Big(k c_{\max} + C_c |\mathcal{S}| |\mathcal{A}| |\T|\Big). \nonumber
    \end{align}
    %We refer to Assumption~\ref{ass: strict feasibility} and Assumption~\ref{ass: continuous} for the definition of the rest symbols.
\end{theorem}

In practice, as optimization error is inevitable, the optimized policy may violate the constraints, \textit{i.e.}, infeasible under the threshold $\gamma_0$. To guarantee the feasibility, one approach is to optimize CMFOMO with $\gamma_0$ replaced by $\gamma_0 - \epsilon_0$ where $\epsilon_0 > 0$ is a small constant. However, the improved feasibility comes at the expense of a worse optimality, as specified in the following theorem.

\begin{theorem} \label{thm: gamma_eps}
    Under Assumption~\ref{ass: strict feasibility} and Assumption~\ref{ass: Lip_continuous}, let $0 < \epsilon_0 < \delta$ be a constant such that Assumption~\ref{ass: strict feasibility} remains satisfied with $\gamma_0$ (\textnormal{resp.} $\delta$) replaced by $\gamma_0 - \epsilon_0$ (\textnormal{resp.} $\delta - \epsilon_0$). Suppose $(L, y, z, \lambda, w)$ is a feasible solution of CMFOMO defined in Theorem~\ref{thm: cmfomo}, with $\gamma_0$ (\textnormal{resp.} $\delta$)  replaced by $\gamma_0 - \epsilon_0$ (\textnormal{resp.} $\delta - \epsilon_0$), and its objective \eqref{eq: cmfomo obj} $\leq \epsilon_0 / \zeta_3$ where $\zeta_3$ is defined by \eqref{eq: zeta_3}. Then, for any $\pi \in \Pi(L)$, its feasibility gap (under $\gamma_0$) is 0, and its optimality gap (under $\gamma_0$) is bounded by:
    \begin{align*}
        0 \leq G_{\text{opt}}(\pi) \leq \frac{\epsilon_0}{\zeta_3} \left[(|\T| + 1) \tilde{\zeta}_1 + \tilde{\zeta}_2 + \tilde{\zeta_4} + \frac{\zeta_3 k|\T|r_{\max}}{\delta - \epsilon_0} \right],
    \end{align*}
    where 
    \begin{align*}
        \tilde{\beta}_y & := \frac{|\mathcal{S}| |\T| (|\T| + 1)}{2} r_{\max} \left(1 + \frac{|\T|}{\delta - \epsilon_0} c_{\max}\right), \\
        \tilde{\beta}_z & := |\mathcal{S}| |\mathcal{A}| (|\T|^2 - |\T| + 2) r_{\max} \left(1 + \frac{|\T|}{\delta - \epsilon_0} c_{\max}\right), \\
        \tilde{\zeta}_1(c_1, c_2) & := \frac{1}{c_1} + \alpha(c_2) \left( C_p \tilde{\beta}_y + C_r + k C_c \frac{|\T| r_{\max}}{\delta - \epsilon_0} \right), \\
        \tilde{\zeta}_2(c_2, c_3) & := \frac{1}{c_3} + \tilde{\beta}_z \alpha(c_2), \\
        \tilde{\zeta}_4(c_2, c_5) & := \frac{1}{c_5} + k \alpha(c_2) \frac{|\T| r_{\max}}{\delta - \epsilon_0} \Big(k c_{\max} + C_c |\mathcal{S}| |\mathcal{A}| |\T|\Big),
    \end{align*}
    and $\alpha(c_2)$ is defined by \eqref{eq: alpha}.
\end{theorem}

The proof of Theorem~\ref{thm: gamma_eps} utilizes similar techniques as the proof of Theorem~\ref{thm: suboptimal cmfomo}, which can also be found in Appendix~\ref{app: proofs}. In addition, as a direct corollary of Theorem~\ref{thm: suboptimal cmfomo}, we get the following asymptotic results.

\begin{corollary}
   Under Assumptions~\ref{ass: strict feasibility} and \ref{ass: Lip_continuous}. Suppose $\left\{\left(c_1^{(i)}, c_2^{(i)}, c_3^{(i)}, c_4^{(i)}, c_5^{(i)}\right)\right\}_{i \geq 1}$ is a sequence of coefficients in CMFOMO's objective function \eqref{eq: cmfomo obj}. Denote by $\left\{\left(L^{(i)}, y^{(i)}, z^{(i)}, \lambda^{(i)}, w^{(i)} \right)\right\}_{i \geq 1}$ a sequence of feasible variables (\textnormal{i.e.}, satisfying all the constraints in Theorem~\ref{thm: cmfomo}) such that:
    \begin{multline*}
        \exists \epsilon > 0, \text{ s.t. } \forall i \geq 1, \text{ } c_1^{(i)} \left|\left| A_{L^{(i)}}^\top y^{(i)} + z^{(i)} + r_{L^{(i)}} - c_{L^{(i)}}^\top \lambda^{(i)}\right|\right|_2 + c_2^{(i)} \left|\left|A_{L^{(i)}} L^{(i)} - b\right|\right|_2 \\ 
        + c_3^{(i)} \left(\left(z^{(i)}\right)^\top L^{(i)}\right) + c_4^{(i)} \left|\left|\gamma_0 - c_{L^{(i)}} L^{(i)} - w^{(i)} \right|\right|_2 + c_5^{(i)} \left|\left(\lambda^{(i)}\right)^\top \left(\gamma_0 - c_{L^{(i)}} L^{(i)}\right)\right| \leq \epsilon.
    \end{multline*}
    We have:
    \begin{enumerate}[1)]
        \item if $\overline{\lim}_{i \rightarrow \infty} c^{(i)}_2 = \overline{\lim}_{i \rightarrow \infty} c^{(i)}_4 = \infty$, then there exists a sub-sequence of $\left\{\left(L^{(i)} \right)\right\}_{i \geq 1}$, whose limiting point is denoted by $L$, such that $\forall \pi \in \Pi(L)$, $\fea(\pi) = 0$.
        \item if $\overline{\lim}_{i \rightarrow \infty} c^{(i)}_1 = \overline{\lim}_{i \rightarrow \infty} c^{(i)}_2 = \overline{\lim}_{i \rightarrow \infty} c^{(i)}_3 = \overline{\lim}_{i \rightarrow \infty} c^{(i)}_4 = \overline{\lim}_{i \rightarrow \infty} c^{(i)}_5 = \infty$, then there exists a sub-sequence of $\left\{\left(L^{(i)} \right)\right\}_{i \geq 1}$, whose limiting point is denoted by $L$, such that $\forall \pi \in \Pi(L)$, $(\pi, L)$ is a CMFG NE defined by Definition~\ref{def: constrained ne}.
    \end{enumerate}
\end{corollary}

Finally, we mention that the inverse of Theorem~\ref{thm: suboptimal cmfomo} is also true, \textit{i.e.}, for any policy $\pi$ whose optimality and feasibility gaps are sufficiently close to 0, one can always find a feasible solution of CMFOMO whose objective is also arbitrarily close to 0. We shall call this property the \textit{characteristicity} of CMFOMO. The proof of the following theorem is also deferred to Appendix~\ref{app: proofs}.

\begin{theorem}[Characteristicity of CMFOMO] \label{thm: char}
    Let $\epsilon_1, \epsilon_2 \geq 0$ and fix $c_i > 0$ $(i = 1, 2, 3, 4, 5)$. Under Assumption~\ref{ass: strict feasibility} and Assumption~\ref{ass: continuous}, let $\pi$ be a policy such that $\left| \opt(\pi) \right| \leq \epsilon_1$ and $0 \leq \fea(\pi) \leq \epsilon_2$. Then, by choosing $L = \Psi(\pi)$, $w = \max\left(0, \gamma_0 - c_L L\right)$ and $(d, y, z, \lambda)$ to be a solution of $\K(L)$ (\textnormal{i.e.}, the \hyperref[kkt]{KKT Conditions} induced by $L$), $(L, y, z, \lambda, w)$ will be a feasible solution of CMFOMO defined in Theorem~\ref{thm: cmfomo} with: 
    \begin{equation*}
        0 \leq \text{objective \eqref{eq: cmfomo obj}} \leq \epsilon_1 (c_3 + 2 c_5) + \epsilon_2 \left[\sqrt{k}\frac{|\T|\rmax}{\delta} (c_3 + c_5) + c_4 \right].
    \end{equation*}
\end{theorem}

\section{Population-Level Constrained Mean-Field Games} \label{sec: pop-level}

In this section, we consider a special class of CMFGs where the cost function is independent of agent state and action, without the strictly feasibility assumption. As we shall see in the next, such CMFGs are degenerate since the representative agent now solves an unconstrained MDP problem. Following a similar procedure as Section~\ref{sec: agent-population-level}, we show that there is a one-to-one correspondence between the CMFG NEs and the optimal solutions of an optimization problem. Finally, we close this section by showing that, roughly speaking, population-level CMFG NEs form a subset of agent-population-level CMFG NEs.

Despite their similarity, we emphasize that Section~\ref{sec: pop-level} is not a special case of Section~\ref{sec: agent-population-level}, as Assumption~\ref{ass: strict feasibility} is no longer imposed. Applying this assumption in the population-level setting would render the constraints redundant, since they are automatically satisfied for any admissible policy. %Because of the same reason, we do not prove the existence of CMFG NEs under population-level constraints.

\subsection{Optimization Framework for Solving NE in CMFGs}

Firstly, under population-level constraints, the CMDP problem $\cmdp{L}$ reduces to an unconstrained MDP problem (provided that the constraints are satisfied under $L$). We denote by $\mdp{L}$ the following unconstrained version of $\cmdp{L}$:
\begin{mdps}
    \begin{align}
        \text{maximize}_{\{ \pi_t \}_{t \in \T}} & \quad \E \left[ \sum_{t \in \T} r_t(s_t, a_t, L_t) \,\bigg|\, s_0 \sim \mu_0 \right] \nonumber \\
        \text{subject to} & \quad a_t \sim \pi_t(s_t), \quad s_{t + 1} \sim p_t(s_t, a_t, L_t). \nonumber
    \end{align}
\end{mdps}
\noindent
Similar to Lemma~\ref{lemma: lp}, $\mdp{L}$ can be transformed into an LP problem (\textit{i.e.}, one only needs to remove the constraints $c_L d \leq \gamma_0$ with the rest remaining the same), and it has the following KKT conditions:
\begin{popkkt} \label{pop_kkt}
    \begin{align*}
        \text{find} \quad & \left\{d \in \R^{|\state||\action||\T|}, y \in \R^{|\state||\T|}, z \in \R^{|\state||\action||\T|} \right\} \nonumber \\
        \text{such that} \quad & -r_L = A_L^\top y + z, \quad A_L d = b, \quad d \geq 0, \quad z \geq 0, \quad z^\top d = 0. 
    \end{align*}
\end{popkkt}
\noindent
In the rest, we denote the above KKT conditions by $\mathcal{K}^p (L)$.

Then, a careful check reveals that, under population-level constraints, $(\pi, L)$ is a CMFG NE defined by Definition~\ref{def: constrained ne} if and only if $(\pi, L)$ is an MFG NE and the constraints are satisfied under $L$:

\begin{proposition} \label{prop: pop cmfomo ne}
    Suppose the constraints are population-level. $(\pi, L)$ is a CMFG NE defined by Definition~\ref{def: constrained ne} if and only if:
    \begin{enumerate}[1)]
        \item $\pi$ is an optimal policy of $\mdp{L}$;
        \item $L = \Psi(\pi)$;
        \item $L$ satisfies \eqref{eq: cmdp3}.
    \end{enumerate}
%     A CMFG NE exists if and only if $\gamma_0$ is greater than or equal to the infimum over all unconstrained MFG NEs. {\color{red} this sentence is incomplete...
%     \[
%     \gamma_0\geq \inf_{(\pi^*,L^*)\in NE} \sum_{t\in \T}c_t(L_t^*),
%     \]
% where $NE$ is the set of all NEs of the unconstrained MFG.
    
%     need to show $\inf$ is achievable? need continuity assumption.
%     }
\end{proposition}

\begin{remark} \label{rmk: pop ne}
     By Proposition~\ref{prop: pop cmfomo ne}, a CMFG NE exists if and only if there exists an MFG NE $(\pi^0, L^0)$ such that
     $\gamma_0 \geq \sum_{t \in \T} c_t(L^0_t)$.
  
     %However, it is difficult to establish a stronger criterion on the existence of CMFG NEs under population-level constraints. In Section~\ref{sec: agent-population-level}, the existence theorem relies on Assumption~\ref{ass: strict feasibility}. But in this context, imposing Assumption~\ref{ass: strict feasibility} would imply that the population-level constraints are always satisfied under any mean-field flow, which makes the CMFGs the same as classical unconstrained MFGs.
\end{remark}

% \begin{remark}
%     Note that if the constraints are population-level, imposing Assumption~\ref{ass: strict feasibility} will make the constraints trivial, \textnormal{i.e.}, the constraints will always be satisfied no matter what value the mean-field flow $L$ takes.
% \end{remark}

% Due to the degenerate nature of population-level constraints, some of the analysis in \S\ref{sec: agent-population-level} is not applicable. For instance, imposing Assumption~\ref{ass: strict feasibility} will imply that the constraints are trivial, \textit{i.e.}, any mean-field flow will satisfy. Alternatively, we propose the following variant of CMFOMO defined in Theorem~\ref{thm: cmfomo} for solving CMFG NEs under population-level constraints.

% Establishing the existence of constrained NE for MFGs with population-level constraints under Assumption~\ref{ass: strict feasibility} is trivial.  Assumption~\ref{ass: strict feasibility} ensures that the population-level constraints will always be satisfied for all mean-field flow $L$. Then Assumption~\ref{ass: continuous} guarantees the existence of NEs and therefore the existence of constrained NEs. Moreover, any NE would be a constrained NE.Therefore throughout  this section, we do not impose Assumption~\ref{ass: strict feasibility}.

Next, we present a result parallel to Theorem~\ref{thm: cmfomo}, which characterizes the CMFG NEs as the optimal solutions to a single optimization problem. We note that Assumption~\ref{ass: strict feasibility} is not imposed in the rest of this section.
\begin{theorem}[Population-level CMFOMO] \label{thm: pop cmfomo}
    Suppose the constraints are population-level. Let $c_i > 0$ $(i = 1, 2, 3, 4)$. $(\pi^\ast, L^\ast)$ is a CMFG NE defined by Definition~\ref{def: constrained ne} if and only if there exists $(y^\ast, z^\ast, w^\ast)$ such that $(L^\ast, y^\ast, z^\ast, w^\ast)$ is an optimal solution of the following optimization problem with the optimal objective value being 0:
    \begin{align}
        \text{minimize}_{L, y, z, w} \quad & c_1 \left|\left| A_L^\top y + z + r_L\right|\right|_2 + c_2 \left|\left|A_L L - b\right|\right|_2 + c_3 (z^\top L) + c_4 \left|\left|\gamma_0 - c_L L - w\right|\right|_2 \label{eq: pop cmfomo obj} \tag{CMFOMO population} \\
        \text{subject to} \quad & L \geq 0, \quad \rv{1}^\top L_t = 1, \quad \forall t \in \T, \nonumber \\
        & ||y||_1 \leq \frac{|\mathcal{S}| |\T| (|\T| + 1)}{2} r_{\max}, \label{eq: pop cmfomo y} \\
        & z \geq 0, \quad ||z||_1 \leq |\mathcal{S}| |\mathcal{A}| \left(|\T|^2 - |\T| + 2\right) r_{\max}, \label{eq: pop cmfomo z} \\
        & 0 \leq w \leq \gamma_0. \nonumber
    \end{align}
\end{theorem}

% \begin{remark}
%     {\color{red}The above optimization problem can be viewed as a special case of \eqref{eq: cmfomo obj} defined in Theorem~\ref{thm: cmfomo} by setting $\lambda = 0$. Intuitively, as $\lambda$ represents the dual variables of the constraints in the LP, this implies that the constraints are redundant in determining the optimality, \textit{i.e.}, any optimal basic solution of the constrained LP is also optimal in the unconstrained LP. See \textit{e.g.} \cite{bertsimas1997introduction} for more details.}
% \end{remark}

\begin{remark}
     Theorem~\ref{thm: pop cmfomo} provides an alternative way to identify the existence of CMFG NEs without Assumption~\ref{ass: strict feasibility}, that is, a CMFG NE exists if and only if the optimal objective value of the above optimization problem is 0. 
\end{remark}

The proof of Theorem~\ref{thm: pop cmfomo} follows directly by combining Proposition~\ref{prop: pop cmfomo ne} and \cite[Theorem 5]{mfomo}. And we mention that the upper bounds on the RHS of \eqref{eq: pop cmfomo y} and \eqref{eq: pop cmfomo z} come from \cite[Corollary 7]{mfomo}. %An application of Theorem~\ref{thm: pop cmfomo} is to identify the existence of CMFG NEs (note that Theorem~\ref{thm: exist} fails without Assumption~\ref{ass: strict feasibility}). 

Let $\pi$ be an arbitrary policy. Its optimality gap (\textit{cf.} Definition~\ref{def: opt_gap}) is only well-defined when $\Psi(\pi)$ satisfies the population-level constraints, \textit{i.e.}, $\fea(\pi) = 0$. Therefore, the bounds \eqref{eq: opt estimate} on $\opt$ in Theorem~\ref{thm: suboptimal cmfomo} are no longer valid. Alternatively, we have the following modified version of Theorem~\ref{thm: suboptimal cmfomo}, whose proof is similar to that of Theorem~\ref{thm: suboptimal cmfomo}.

\begin{theorem}[Suboptimality of population-level CMFOMO] 
    Let $\epsilon \geq 0$. Under Assumption~\ref{ass: continuous}, suppose $(L, y, z, w)$ is a feasible solution of CMFOMO defined in Theorem~\ref{thm: pop cmfomo} with its objective 
    \[
    0 \leq \eqref{eq: pop cmfomo obj} \leq \epsilon.
    \]
    Then, for any $\pi \in \Pi(L)$, its feasibility gap is bounded by:
    \begin{align*}
        0 \leq G_{\text{fea}}(\pi) & \leq \epsilon \zeta_3,
    \end{align*}
    where $\zeta_3$ is defined in \eqref{eq: zeta_3}. If specially $G_{\text{fea}}(\pi) = 0$, then the optimality gap will be well-defined and bounded by:
    \begin{equation*}
        0 \leq G_{\text{opt}}(\pi) \leq \epsilon \Big[(|\T| + 1) \zeta_1^\prime + \zeta_2^\prime \Big],
    \end{equation*}
    where
    \begin{align*}
        \beta_y^\prime & := \frac{|\mathcal{S}| |\T| (|\T| + 1)}{2} r_{\max}, \\
        \beta_z^\prime & := |\mathcal{S}| |\mathcal{A}| (|\T|^2 - |\T| + 2) r_{\max}, \\
        \zeta_1^\prime(c_1, c_2) & := \frac{1}{c_1} + \alpha(c_2) \left( C_p \beta_y^\prime + C_r \right), \\
        \zeta_2^\prime(c_2, c_3) & := \frac{1}{c_3} + \beta_z^\prime \alpha(c_2),
    \end{align*}
    and $\alpha$ is defined in \eqref{eq: alpha}.
\end{theorem}

A similar result to Theorem~\ref{thm: char} can also be established by requiring $\epsilon_2 = 0$, as formally stated below.

\begin{theorem}[Characteristicity under population-level constraints]
    Suppose the constraints are population-level. Let $\epsilon_1 \geq 0$ and fix $c_i > 0$ $(i = 1, 2, 3, 4)$. Let $\pi$ be a policy such that $\fea(\pi) = 0$ and $0 \leq \opt(\pi) \leq \epsilon_1$. Then, by choosing $L = \Psi(\pi)$, $w = 0$ and $(d, y, z)$ to be a solution of $\K^p(L)$ (\textnormal{i.e.}, the \hyperref[pop_kkt]{KKT Conditions} induced by $L$), $(L, y, z, w)$ will be a feasible solution of CMFOMO defined in Theorem~\ref{thm: pop cmfomo} with its objective: 
    \begin{equation*}
        0 \leq \eqref{eq: pop cmfomo obj} \leq \epsilon_1 c_3.
    \end{equation*}
\end{theorem}

\subsection{Connection between Population-Level and Agent-Population-Level Constraints}

We end this section with a discussion on the connection between population-level and agent-population-level constraints. Our main result is that population-level CMFG NEs comprise of a subset of agent-population-level CMFG NEs.

\begin{definition} \label{def: twinned}
    Let $c^p$ (\textnormal{resp.} $c^a$) be the cost function of population-level (\textnormal{resp.} agent-population-level) constraints. $c^p$ and $c^a$ are said to be \textnormal{twinned} if for any mean-field flow $L$, it holds that:
    \[
    \sum_{t \in \T} c^p_t\big(L_t\big)=\sum_{t \in \T} \sum_{s \in \state, a \in \action} c^a_t\big(s,a,L_t\big)L_t(s,a).
    \]
    % \begin{equation*}
    %     \sum_{t \in \T} c^p_t\big(\Psi_t(\pi)\big) = \E_{a \sim \pi} \left[ \sum_{t \in \T} c^a_t\big(s_t, a_t, \Psi_t(\pi)\big) \,\bigg|\, s_0 \sim \mu_0 \right].
    % \end{equation*}
    % {\color{red} if for any $L$: $\sum_{t \in \T} c^p_t\big(L_t\big)=\sum_{t \in \T} \sum_{s,a}c^a_t\big(s,a,L_t\big)L_t(s,a)$?}
\end{definition}

The constraints \eqref{eq: example_pop} and \eqref{eq: example_agent} in Example~\ref{exam: SIS} are an example of twinned cost functions. The point of Definition~\ref{def: twinned} is that for any policy $\pi$, if $\pi$ is feasible in $\mathcal{M}_{c^p}\big(\Psi(\pi)\big)$, then $\pi$ will also be feasible in $\mathcal{M}_{c^a}\big(\Psi(\pi)\big)$. Furthermore, if $\pi$ is optimal in $\mathcal{M}_{c^p}\big(\Psi(\pi)\big)$ (which is equivalent to $\pi$ being optimal in $\mathcal{M}\big(\Psi(\pi)\big)$), then $\pi$ will also be optimal in $\mathcal{M}_{c^a}\big(\Psi(\pi)\big)$. This argument proves the following Theorem~\ref{thm: subset}. And we mention that the calculation results in Example~\ref{exam: SIS} conform to Theorem~\ref{thm: subset}.

\begin{theorem} \label{thm: subset}
    Let $c^p$ and $c^a$ be a pair of twinned cost functions at the population level and agent population level, respectively. If $(\pi, L)$ is a CMFG NE defined by Definition~\ref{def: constrained ne} under $c^p$, then $(\pi, L)$ will also be a CMFG NE under $c^a$.
\end{theorem}

\section{Approximation of NE in Finite-player Games with Agent-population-level Constraints} \label{sec: N_player}

It is well-known that NEs of MFGs are approximated NEs of $N$-player games. In this section, we provide the corresponding analysis and show that NEs of CMFGs are also approximated NEs of the counterpart constrained $N$-player games, with an explicit dependency of the approximation error on the number of players $N$. Throughout this section, we always assume that the constraints are agent-population-level (\textit{cf.} Section~\ref{sec: agent-population-level}).

Let $0 < N < +\infty$.
%In this section, we investigate the NE in games where there are $N$ players. 
We assume that the $N$ players are homogeneous and have mean-field-type interactions. %, and 2) the $N$ players are restricted to only take policies that satisfy certain constraints. The main point of this section is that CMFGs can be viewed as an approximation to $N$-player games when $N$ tends to infinity. 
In the following, we prove that the NE of CMFGs (\textit{cf.} Definition~\ref{def: constrained ne}) is an $\epsilon$-NE of $N$-player games, and that $\epsilon$ can be made arbitrarily close to 0 as $N$ tends to infinity. In the rest, we inherit the notations in \S\ref{sec: setup}, and we use the superscript $^{(i)}$ to indicate the private variables of the $i$-th player.

To begin with, for each $1 \leq i \leq N$, the constrained control problem of the $i$-th player writes:

\begin{N_cmdps} \label{problem: N_cmdps}
    \begin{align}
        \text{maximize}_{\left\{ \pi^{(i)}_t \right\}_{t \in \T}} & \quad \E \left[ \sum_{t \in \T} r_t\left(s^{(i)}_t, a^{(i)}_t, L^{(N)}_t\right) \,\bigg|\, s^{(i)}_0 \sim \mu_0 \right] \label{eq: N_cmdp1} \tag{$N$ reward} \\
        \text{subject to} & \quad a^{(i)}_t \sim \pi^{(i)}_t\left(s^{(i)}_t\right), \quad s^{(i)}_{t + 1} \sim p_t\left(s^{(i)}_t, a^{(i)}_t, L^{(N)}_t\right), \tag{$N$ dynamics} \\
        & \quad \E \left[ \sum_{t \in \T} c_t\left(s^{(i)}_t, a^{(i)}_t, L^{(N)}_t\right) \,\bigg|\, s^{(i)}_0 \sim \mu_0 \right] \leq \gamma_0, \tag{$N$ constraints}
    \end{align}
\end{N_cmdps}
\noindent
where $L^{(N)}_t$ denotes the empirical distribution of the $N$ players at time step $t$, \textit{i.e.},
\begin{equation*}
    \forall t \in \T, s \in \state, a \in \action, \quad L^{(N)}_t (s, a) := \frac{1}{N} \sum_{i=1}^N \delta_{\left(s^{(i)}_t, a^{(i)}_t\right) = (s, a)}.
\end{equation*}
For convenience, with an abuse of notation, we use $\cmdpN{i}$ to denote the above constrained control problem of the $i$-th player, and we denote by $V^{(i)}\left(\pi^{(1)}, \cdots, \pi^{(i)}, \cdots, \pi^{(N)}\right)$ the expected cumulative reward of the $i$-th player when the players take policies $\left\{\pi^{(i)}\right\}_{1 \leq i \leq N}$.

\begin{remark}
    Unlike that the players are independent in the constrained MDP problem $\cmdp{L}$ defined in \S\ref{sec: setup}, the $N$ players in $\left\{\cmdpN{i}\right\}_{1 \leq i \leq N}$ are correlated. This is because the decision of a specific player does affect the law of the empirical distribution of the $N$ players, which in turn affects the evolution of other players. Another consequence is that $\cmdpN{i}$ is no longer a traditional MDP even when the policies of the rest players are fixed. 
\end{remark}

Next, we generalize the concept of feasibility gap (\textit{cf.} Definition~\ref{def: fea_gap}) to the $N$-player setting.

\begin{definition}[Feasibility gap of the $i$-th player] \label{def: fea_gap_N}
    Denote by $\left\{\pi^{(i)}\right\}_{1 \leq i \leq N}$ the policies of the $N$ players. For any $1 \leq i \leq N$, we define the feasibility gap of the $i$-th player by:
    \begin{equation*}
        G_{\text{fea}}^{(i)}\left(\pi^{(1)}, \cdots, \pi^{(i)}, \cdots, \pi^{(N)}\right) := \left|\left|\min\Big(0, \gamma_0 - \gamma^{(i)}\left(\pi^{(1)}, \cdots, \pi^{(i)}, \cdots, \pi^{(N)}\right) \Big)\right|\right|_2,
    \end{equation*}
    where $\min$ is taken entry-wise, and
    \begin{equation*}
        \gamma^{(i)}\left(\pi^{(1)}, \cdots, \pi^{(i)}, \cdots, \pi^{(N)}\right) := \E_{\forall 1 \leq j \leq N,\, a^{(j)} \sim \pi^{(j)}} \left[ \sum_{t \in \T} c_t\left(s^{(i)}_t, a^{(i)}_t, L^{(N)}_t\right) \,\bigg|\, s^{(i)}_0 \sim \mu_0 \right].
    \end{equation*}
\end{definition}

Given the policies of all but the $i$-th player, we denote by $\F^{(i)}\left(\pi^{(1)}, \cdots, \pi^{(i-1)}, \pi^{(i+1)}, \cdots, \pi^{(N)}\right)$ the set of policies of the $i$-th player whose feasibility gap is 0 (\textit{i.e.}, the set of feasible policies of the $i$-th player provided the policies of the rest players). To ease the exposure of notations, we abbreviate it by $\F^{(i)}$ if the policies of the rest players are clearly indicated in the context, similarly for $G_{\text{fea}}^{(i)}$. 

\begin{definition}[$\epsilon$-NE of $N$-player games] \label{def: eps-NE}
    Let $\epsilon_1, \epsilon_2 \geq 0$. The set of policies $\left\{\pi^{(i)}\right\}_{1 \leq i \leq N}$ is called an $(\epsilon_1, \epsilon_2)$-NE of the $N$-player game $\left\{\cmdpN{i}\right\}_{1 \leq i \leq N}$ if for any $1 \leq i \leq N$:
    \begin{align*}
        \sup_{\pi \in \F^{(i)}} V^{(i)}\left(\pi^{(1)}, \cdots, \pi^{(i-1)}, \pi, \pi^{(i+1)}, \cdots, \pi^{(N)}\right)
        \leq V^{(i)}\left(\pi^{(1)}, \cdots, \pi^{(i-1)}, \pi^{(i)}, \pi^{(i+1)}, \cdots, \pi^{(N)}\right) + \epsilon_1.
    \end{align*}
    and $G_{\text{fea}}^{(i)} \leq \epsilon_2$. Specially, $\left\{\pi^{(i)}\right\}_{1 \leq i \leq N}$ is called an NE if $\epsilon_1 = \epsilon_2 = 0$.
\end{definition}

However, Definition~\ref{def: eps-NE} is not well-defined when $\F^{(i)}$ is empty. To ensure that $\F^{(i)}$ is non-empty (especially when $N$ is large enough), we need to slightly strengthen Assumption~\ref{ass: strict feasibility} by the following assumption. 

\begin{assumption}[Strengthened strict feasibility] \label{ass: strict feasibility N}
    There exists $\delta > 0$, such that for any mean-field flow $L$, there exists a policy $\pi$ that satisfies:
    \begin{equation*}
        \E_{a \sim \pi} \left[ \sum_{t \in \T} c_t (s_t, a_t, L_t) \,\Big|\, s_0 \sim \mu_0 \right] \leq \gamma_0 - \delta.
    \end{equation*}
\end{assumption}

\begin{remark}
    Assumption~\ref{ass: strict feasibility} only requires the existence of strictly feasible policy for each of the constraints. The feasible policies can be different across different constraints. However, in Assumption~\ref{ass: strict feasibility N}, we require the existence of a single  policy that is strictly feasible for all constraints. When there is only one constraint, \textit{i.e.}, when $k=1$, the two assumptions are equivalent.
\end{remark}

With the new assumption, we present the main result of this section, whose proof can be found in Appendix~\ref{app: N_player}.

\begin{theorem}[$\epsilon$-NE approximated by CMFGs] \label{thm: eps-NE}
    Under Assumptions~\ref{ass: Lip_continuous} and \ref{ass: strict feasibility N}, suppose $(\pi^\ast, L^\ast)$ is a CMFG NE defined by Definition~\ref{def: constrained ne}. Fix $0 < \epsilon \leq \delta$ and let $N > 0$ be such that:
    \begin{equation*}
        \left(\frac{1}{2\sqrt{N}} + \frac{2}{N}\right) \tilde{C} \leq \epsilon,
    \end{equation*}
    where
    \begin{align*}
         \tilde{C} & := C_p |\state| |\action| (|\T| - 1)^2 \frac{C_{p, \state, \action}^{|\T|} - 1}{C_{p, \state, \action} - 1} \max\Big\{C_r + r_{\max},\, C_c + c_{\max}\Big\}, \\
         C_{p, \state, \action} & := \left(C_p + 1\right) |\state| |\action|.
    \end{align*}
    Then, $\left\{\pi^\ast\right\}_{1 \leq i \leq N}$ (\textnormal{i.e.}, $N$ copies of $\pi^\ast$) is an $(\epsilon_1, \epsilon_2)$-NE of the $N$-player game defined by Definition~\ref{def: eps-NE}, where
    \begin{align*}
        \epsilon_1 & = \epsilon \left(2 + \frac{k |\T| r_{\max}}{\delta}\right), \\
        \epsilon_2 & = \epsilon \sqrt{k}.
    \end{align*}
\end{theorem}

Let $\epsilon_0 > 0$ be a small constant. As a corollary of Theorem~\ref{thm: eps-NE}, any CMFG NE under $\gamma_0 - \epsilon_0$ will be strictly feasible in the finite-player game under $\gamma_0$ when the number of players is large enough. Same as Theorem~\ref{thm: gamma_eps}, the improved feasibility comes at the expense of a worse optimality.

\begin{corollary}
    Under Assumptions~\ref{ass: Lip_continuous} and \ref{ass: strict feasibility N}, let $0 < \epsilon_0 < \delta$ be a constant. Suppose $(\pi^\ast, L^\ast)$ is a CMFG NE defined by Definition~\ref{def: constrained ne} by replacing $\gamma_0$ with $\gamma_0 - \epsilon_0$. Fix $0 < \epsilon \leq \epsilon_0$ and let $N > 0$ be such that:
    \begin{equation*}
        \left(\frac{1}{2\sqrt{N}} + \frac{2}{N}\right) \tilde{C} \leq \epsilon,
    \end{equation*}
    where $\tilde{C}$ is defined in Theorem~\ref{thm: eps-NE}. Then, $\left\{\pi^\ast\right\}_{1 \leq i \leq N}$ (\textnormal{i.e.}, $N$ copies of $\pi^\ast$) is an $(\tilde{\epsilon}_1, 0)$-NE of the $N$-player game (under $\gamma_0$) defined by Definition~\ref{def: eps-NE}, where
    \begin{align*}
        \tilde{\epsilon}_1 = \epsilon_0 \frac{k|\T|r_{\max}}{\delta - \epsilon_0} + \epsilon \left(2 + \frac{k |\T| r_{\max}}{\delta}\right).
    \end{align*}
\end{corollary}

Finally, by connecting Theorem~\ref{thm: eps-NE} with Theorem~\ref{thm: suboptimal cmfomo}, we show that CMFOMO can be used to obtain $\epsilon$-NEs of the finite-player game.

\begin{corollary}[$\epsilon$-NE approximated by CMFOMO]
    Under Assumptions~\ref{ass: Lip_continuous} and \ref{ass: strict feasibility N}, let $\epsilon^\prime > 0$ and suppose $(L, y, z, \lambda, w)$ is a feasible solution of CMFOMO defined in Theorem~\ref{thm: cmfomo} with its objective \eqref{eq: cmfomo obj} $\leq \epsilon^\prime$. Fix $0 < \epsilon \leq \delta$ and let $N$ be such that:
    \begin{equation*}
        \left(\frac{1}{2\sqrt{N}} + \frac{2}{N}\right) \tilde{C} \leq \epsilon,
    \end{equation*}
    where $\tilde{C}$ is defined in Theorem~\ref{thm: eps-NE}. Then, for any $\pi \in \Pi(L)$, one can show that $\left\{\pi\right\}_{1 \leq i \leq N}$ (\textnormal{i.e.}, $N$ copies of $\pi$) is an $(\epsilon^\prime_1, \epsilon^\prime_2)$-NE of the $N$-player game defined by Definition~\ref{def: eps-NE}, where
    \begin{align*}
        \epsilon_1^\prime & :=  \epsilon^\prime \Big[(|\T| + 1) \zeta_1 + \zeta_2 + \zeta_4 \Big] + \epsilon \left(2 + \frac{k |\T| r_{\max}}{\delta}\right), \\
        \epsilon_2^\prime & := \epsilon^\prime \zeta_3 + \epsilon \sqrt{K}.
    \end{align*}
    We refer to Theorem~\ref{thm: suboptimal cmfomo} for the definition of $\zeta_1, \zeta_2, \zeta_3, \zeta_4$.
\end{corollary}

\section{Numerical Experiments} \label{sec: numerical}

In this section, we illustrate our CMFOMO framework developed before on the modified SIS problem introduced in Example~\ref{exam: SIS}. We will discuss different types of constraints including both agent-population-level constraints and population-level constraints and how they affect the solutions of the MFGs. Our implementation is based on PyTorch and MFGLib proposed in \cite{mfglib}. %games modified from the literature. Various micro and macro constraints are investigated (\cf Remark~\ref{rmk: macro_micro}).

%\subsection*{Constrained SIS model (\emph{cf.} Example~\ref{exam: SIS}).}

In the following numerical experiments, we set $T = 10$ and initialize the population with $\mu_0(s_0 = I) = 0.5$. The reward functions and transition probabilities follow those in Example~\ref{exam: SIS}. Note that, because of the different choices of $T$ and $\mu_0$, the threshold value of $\gamma_0$ ensuring the existence of a Nash equilibrium differs from that in Example~\ref{exam: SIS}.

\subsection{Agent-population-level constraint on state} \label{sec: agent_state}

In this section, we impose a condition on the agents by restricting the probability of each agent getting infected to be upper bounded by a certain threshold. Mathematically speaking, we consider the following agent-population-level constraint on state:
\begin{equation}
    \frac{1}{|\T|} \E \left( \sum_{t \in \T} \delta_{s_t = I} \,\bigg|\, s_0 \sim \mu_0 \right) \leq \gamma_0. \label{eq: gamma_0.25_con}
\end{equation}
%where we set $\gamma_0 = 0.25$. 
Note that this setting satisfies Assumptions~\ref{ass: strict feasibility} and \ref{ass: continuous}, so there exists a CMFG NE defined by Definition~\ref{def: constrained ne} according to Theorem~\ref{thm: exist}.

We use Projected Gradient Descent (PGD) to optimize the objective function \eqref{eq: cmfomo obj} of CMFOMO defined in Theorem~\ref{thm: cmfomo}, combined with the Adam optimizer embedded in PyTorch. Inspired by Theorem~\ref{thm: char}, in each iteration, once $L$ and $\lambda$ are determined, we solve $y$ and $z$ from the following modified KKT conditions ($d$ is discarded):
\begin{align}
    \text{find} \quad & \left\{d \in \R^{|\state||\action||\T|}, y \in \R^{|\state||\T|}, z \in \R^{|\state||\action||\T|}\right\} \nonumber \\
    \text{such that} \quad & -r_L + c_L^\top\lambda = A_L^\top y + z, \quad A_L d = b, \quad d \geq 0, \quad z \geq 0, \quad z^\top d = 0, \nonumber
\end{align}
and determine $w$ the same way as Theorem~\ref{thm: char}, \textit{i.e.}, $w = \max\left(0, \gamma_0 - c_L L\right)$. This helps our algorithm to converge faster. In a word, the only variables to be optimized are $L$ and $\lambda$, and our optimization procedure is summarized below:
\begin{algorithm}[H]
    \begin{algorithmic}[1]
        \Require Threshold $\gamma_0$; initial values for $(L^0, \lambda^0)$; loss weights $c_i$ $(i = 1, 2, 3, 4, 5)$ for CMFOMO; gradient optimizer; maximum number of iterations $N$.
        \State $n \gets 0$
        \While{$n < N$}
            \State Compute $(y^n, z^n, w^n)$ by the modified KKT conditions and Theorem~\ref{thm: char}
            \State $\text{loss}^n \gets \text{CMFOMO} (L^n, y^n, z^n, \lambda^n, w^n)$ \Comment{Objective \eqref{eq: cmfomo obj} of CMFOMO}
            \State $\text{loss}^n$.backward()
            \State $\hat{L}^{n + 1} \gets \text{optimizer}(n, L^n, L^n\text{.grad}), \quad \hat{\lambda}^{n + 1} \gets \text{optimizer}(n, \lambda^n, \lambda^n\text{.grad})$
            \State $L^{n + 1} \gets \text{proj}_L(\hat{L}^{n + 1}), \quad \lambda^{n + 1} \gets \text{proj}_\lambda(\hat{\lambda}^{n + 1})$ \Comment{Ensure $L^{n+1}$ a mean-field flow and $\lambda^{n+1} \geq 0$} \label{line: proj}
            \State $n \gets n + 1$
        \EndWhile
        \State \Return $L^N$
    \end{algorithmic}
    \caption{Optimization of CMFOMO defined in Theorem~\ref{thm: cmfomo}}
    \label{alg: SIS_opt}
\end{algorithm}

As for other hyper-parameters, the initial value $L^0$ is induced by the uniform policy, \textit{i.e.}, the policy that chooses all actions with an equal probability. The initial value $\lambda^0$ is chosen to be 0. We choose an equal weight for CMFOMO, \textit{i.e.}, $c_1 = c_2 = c_3 = c_4 = c_5 = 1$. Finally, we set the learning rate of the Adam optimizer to be $5 \times 10^{-3}$.

\begin{figure}[H]
    \centering
    \begin{minipage}[b]{0.49\textwidth}
        \includegraphics[width=\textwidth]{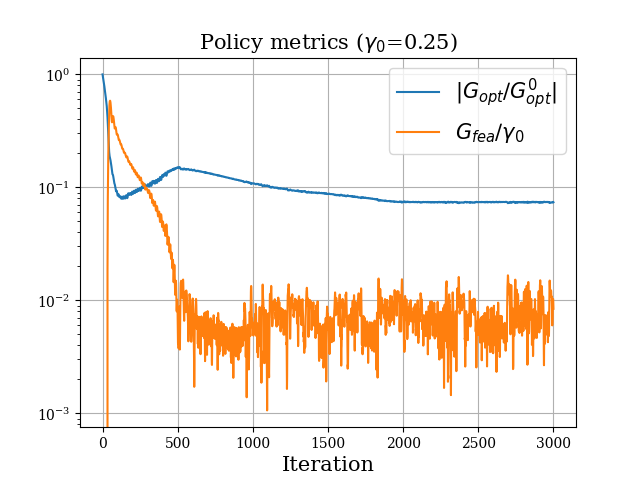}
    \end{minipage}
    \begin{minipage}[b]{0.49\textwidth}
        \includegraphics[width=\textwidth]{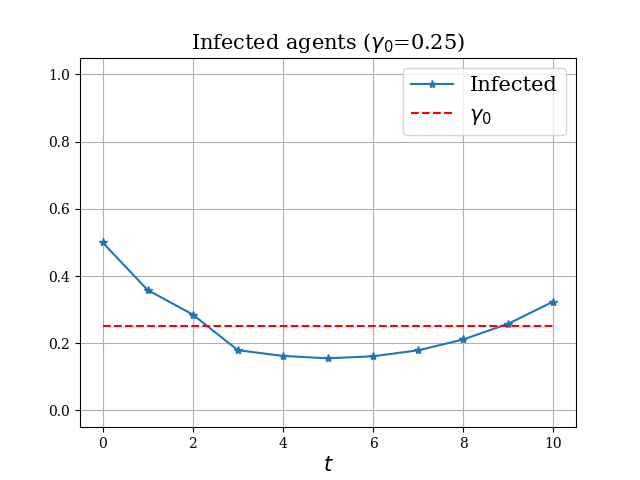}
    \end{minipage}
    \caption{Optimized policy of the constrained SIS model when $\gamma_0 = 0.25$. Left: the evolution of $\opt$ and $\fea$ with respect to optimization iteration, where $\opt^0$ represents the optimality gap of the initial policy; Right: the percentage of infected agents at each time step under the optimized policy, where the red dotted line represents the threshold $\gamma_0$ of the agent-population-level constraint \eqref{eq: gamma_0.25_con}.}
    \label{fig: SIS_gamma_0.25}
\end{figure}

When $\gamma_0$ is set 0.25, in Figure~\ref{fig: SIS_gamma_0.25} (Left), we plot the evolution of $\opt$ and $\fea$ with respect to optimization iteration. In the figure legends, ``$\opt$'' and ``$\fea$'' represent the optimality gap and feasibility gap of the policy at the current iteration, respectively, and ``$\opt^0$'' refers to the optimality gap of the initial policy.  We note that $\fea / \gamma_0$ (see the orange line) is initially very small because the starting policy is feasible. As the iterations proceed, the policy quickly moves toward minimizing $\opt$, temporarily sacrificing feasibility. Nevertheless, $\fea$ eventually converges to an acceptable region. In Figure~\ref{fig: SIS_gamma_0.25} (Right), we plot the values of $\big\{\E\left( \delta_{s_t = I} \right)\big\}_{t \in \T}$ under the optimized policy, \textit{i.e.}, the percentage of infected agents at each time step. Under the optimized policy, the LHS of \eqref{eq: gamma_0.25_con} is $0.2521$, slightly exceeding the threshold $\gamma_0 = 0.25$.

\begin{remark} \label{rmk: eps_0}
    Here the optimized policy is not feasible under $\gamma_0 = 0.25$. By Theorem~\ref{thm: gamma_eps}, in order to produce a feasible policy under $\gamma_0 = 0.25$, one practical way is to optimize CMFOMO with $\gamma_0$ replaced by $\gamma_0 - \epsilon_0$ for some small $\epsilon_0 > 0$. For illustration, in Table~\ref{tab: eps_0}, we vary $\epsilon_0 \in [0.01, 0.05]$ and list the optimality gap (under $\gamma_0 = 0.25$) and average infected agents (\textnormal{i.e.}, the LHS of \eqref{eq: gamma_0.25_con}) of the optimized policies. As observed, with the increasing of $\epsilon_0$, the feasibility is improved at the expense of a worse optimality gap.
\end{remark}

\begin{table}[H]
    \centering
    \begin{tabular}{ccccc}
    \hline
    \hline
    $\epsilon_0$ && $\opt(\pi)$ under $\gamma_0 = 0.25$ && Average infected agents \\
    \hline
    0    && 0.1625 && 0.2521 \\
    0.01 && 0.1642 && 0.2513 \\
    0.02 && 0.1663 && 0.2493 \\
    0.03 && 0.1664 && 0.2486 \\
    0.04 && 0.1675 && 0.2473 \\
    0.05 && 0.1880 && 0.2469 \\
    \hline
    \hline
    \end{tabular}
    \caption{Optimality gap (under $\gamma_0 = 0.25$) and average infected agents (\textit{i.e.}, the LHS of \eqref{eq: gamma_0.25_con}) of the policy optimized under the threshold $\gamma_0 - \epsilon_0$ (\textit{cf.} Remark~\ref{rmk: eps_0}).}
    \label{tab: eps_0}
\end{table}

% Then, we vary the threshold $\gamma_0$ and investigate how the solution changes with different values of $\gamma_0$. In Figure~\ref{fig: SIS_sensitivity} (Left), we plot the percentage of infected agents with the threshold $\gamma_0$ taking different values. In Figure~\ref{fig: SIS_sensitivity} (Right), we plot the percentage of agents who go out (\textit{i.e.} $a_t = U$). We observe that the smaller the threshold $\gamma_0$ is, the more conservative the agents will be, \textit{i.e.} going out less, which then leads to a lower infection rate. In addition, when $\gamma_0$ exceeds a certain value, the optimized policies coincide (\textit{e.g.}, the trajectories corresponding to $\gamma_0 = 0.5$ and 1 overlap well in Figure~\ref{fig: SIS_sensitivity}). 

\begin{figure}[H]
    \centering
    \begin{minipage}[b]{0.49\textwidth}
        \includegraphics[width=\textwidth]{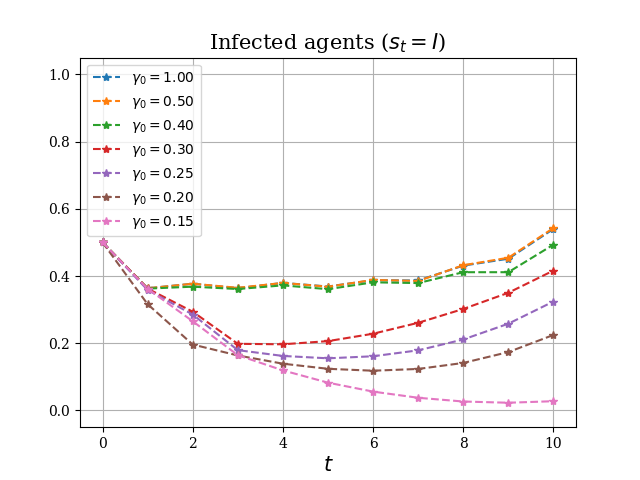}
    \end{minipage}
    \begin{minipage}[b]{0.49\textwidth}
        \includegraphics[width=\textwidth]{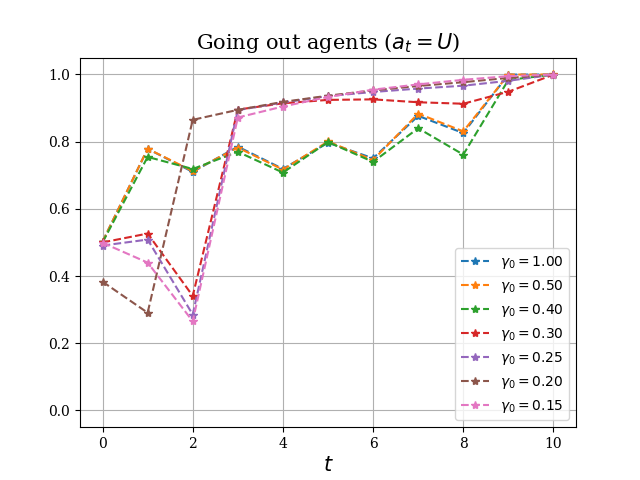}
    \end{minipage}
    \caption{Sensitivity analysis under the agent-population-level constraint \eqref{eq: gamma_0.25_con}. The $x$-axis represents the time steps and the $y$-axis represents the percentage of infected agents (Left) or agents who go out (Right). The policies are optimized when $\gamma_0$ takes different values.}
    \label{fig: SIS_sensitivity}
\end{figure}

We then vary the threshold $\gamma_0$ to examine how the equilibrium changes with different constraint levels. Figure~\ref{fig: SIS_sensitivity} (Left) plots the percentage of infected agents (\textit{i.e.}, $s_t = I$). As expected, the infection rate decreases as $\gamma_0$ becomes smaller, reflecting a stronger constraint. Figure~\ref{fig: SIS_sensitivity} (Right) shows the proportion of agents choosing to go out (\textit{i.e.}, $a_t = U$). Two clear patterns emerge: (1) for small $\gamma_0$ values (e.g., $\gamma_0 \le 0.3$), agents initially prefer to stay home ($a_t = D$) to reduce infection risk, and once the infection rate drops sufficiently, they gradually resume going out to increase rewards; (2) for larger $\gamma_0$, where the constraint is relatively loose, agents go out from the start to maximize rewards. Moreover, we observe that the optimized policies coincide once $\gamma_0$ exceeds a certain threshold; for instance, the trajectories for $\gamma_0 = 0.5$ and $\gamma_0 = 1$ nearly overlap in Figure~\ref{fig: SIS_sensitivity}. This is because the constraint becomes inactive and the solution of the CMFG is the same as the the unconstrained MFG.

\subsection{Agent-population-level constraint on action}

\begin{figure}[H]
    \centering
    \begin{minipage}[b]{0.49\textwidth}
        \includegraphics[width=\textwidth]{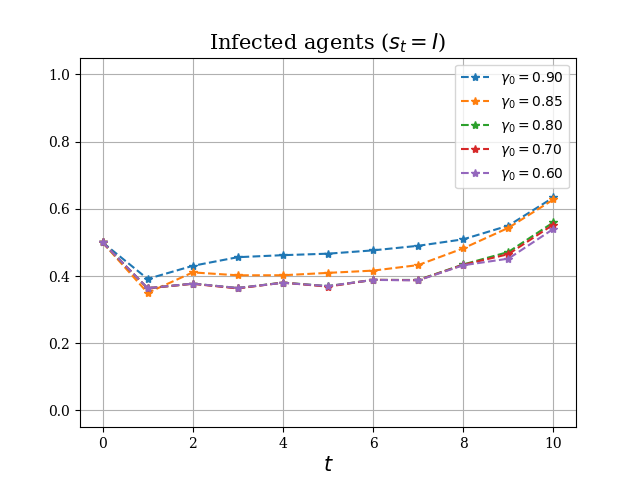}
    \end{minipage}
    \begin{minipage}[b]{0.49\textwidth}
        \includegraphics[width=\textwidth]{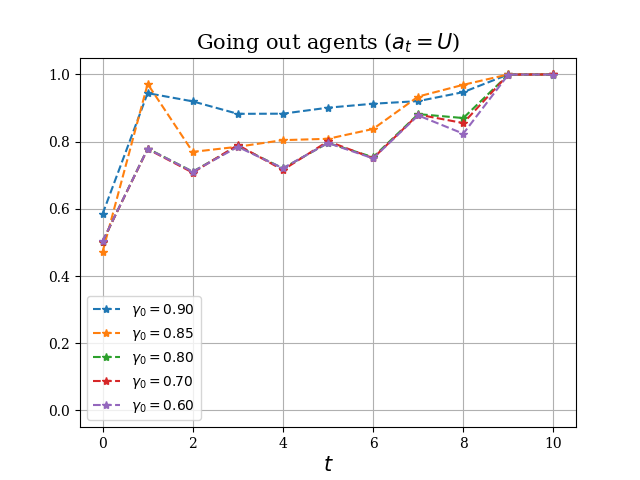}
    \end{minipage}
    \caption{Sensitivity analysis under the agent-population-level constraint \eqref{eq: agent_action}. The $x$-axis represents the time steps and the $y$-axis represents the percentage of infected agents (Left) or agents who go out (Right). The policies are optimized when $\gamma_0$ takes different values.}
    \label{fig: SIS_sensitivity_action}
\end{figure}

In this section, we study a different agent-population-level constraint which guarantees that the probability of each agent going out is lower bounded by a certain threshold.
Specifically, we impose the following agent-population-level constraint on action:
\begin{equation} \label{eq: agent_action}
    \frac{1}{|\T|} \E \left( \sum_{t \in \T} \delta_{a_t = U} \,\bigg|\, s_0 \sim \mu_0 \right) \geq \gamma_0.
\end{equation}
According to Theorem~\ref{thm: exist}, a CMFG NE exist as long as $\gamma_0 < 1$. In Figure~\ref{fig: SIS_sensitivity_action}, we plot the percentage of infected (Left) and going out (Right) agents under different $\gamma_0$ values. The optimization procedure is exactly the same as Section~\ref{sec: agent_state}. As expected, the larger the value of $\gamma_0$, the higher the infection rate will be, \textit{i.e.}, the agents will have to go out more frequently. In addition, when $\gamma_0$ is smaller than a certain value, the optimized policies coincide (\textit{e.g.}, the trajectories corresponding to $\gamma_0 = 0.7$ and 0.6 overlap well in Figure~\ref{fig: SIS_sensitivity_action}). 

\begin{figure}[H]
    \centering
    \includegraphics[width=0.5\textwidth]{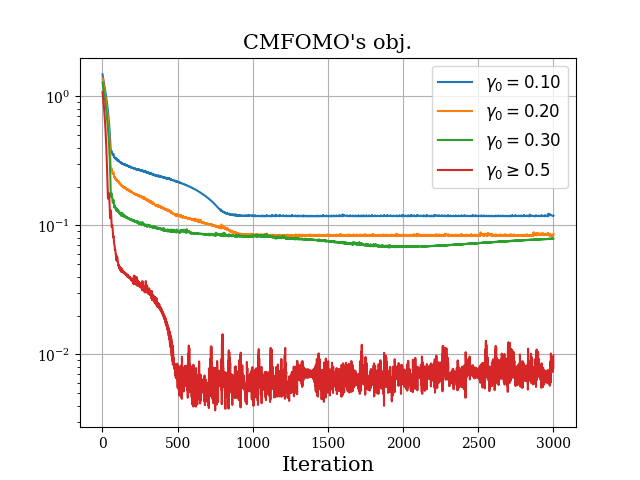}
    \caption{Evolution of the objective function \eqref{eq: pop cmfomo obj} of CMFOMO defined in Theorem~\ref{thm: pop cmfomo} under the population-level constraint \eqref{eq: macro_exp_con}. The $x$-axis represents the optimization iteration, and the $y$-axis represents the objective value. Different lines correspond to different values of $\gamma_0$.}
    \label{fig: macro_loss}
\end{figure}

\subsection{Population-level constraint on state}

Finally, in this section, we consider the following population-level constraint on agent state:
\begin{equation}
    \frac{1}{|\T|} \sum_{t \in \T} L(s_t = I) \leq \gamma_0. \label{eq: macro_exp_con}
\end{equation}
When optimizing the objective function \eqref{eq: pop cmfomo obj} of CMFOMO defined in Theorem~\ref{thm: pop cmfomo}, we devise an algorithm similar to Algorithm~\ref{alg: SIS_opt}. The core idea is that \eqref{eq: pop cmfomo obj} can be viewed as plugging in $\lambda 
= 0$ into \eqref{eq: cmfomo obj}. By Proposition~\ref{prop: pop cmfomo ne} and Remark~\ref{rmk: pop ne}, the objective function \eqref{eq: pop cmfomo obj} admits a zero point if and only if $\gamma_0$ is no smaller than a certain value. 

In Figure~\ref{fig: macro_loss}, we plot the evolution of the value of the objective function \eqref{eq: pop cmfomo obj} with respect to optimization iteration, with $\gamma_0$ taking different values. As shown, when $\gamma_0$ is large enough such that a CMFG NE exists (\textit{e.g.}, $\gamma_0 \geq 0.5$), the objective decreases fast. In contrast, when $\gamma_0$ is smaller such that a CMFG NE does not exist, the objective stays at a relatively higher level.

\begin{figure}[H]
    \centering
    \begin{minipage}[b]{0.49\textwidth}
        \includegraphics[width=\textwidth]{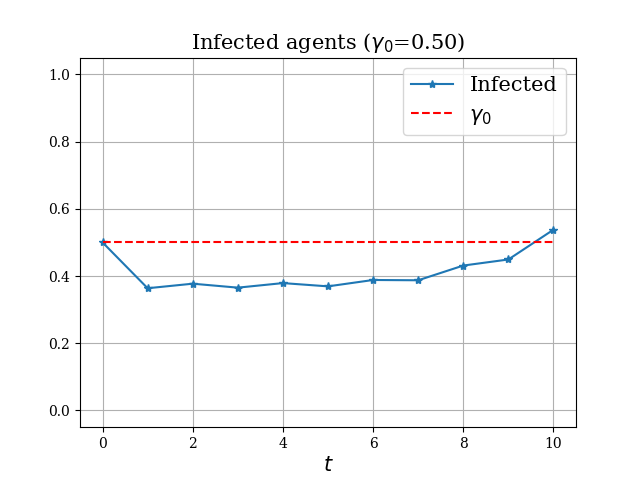}
    \end{minipage}
    \begin{minipage}[b]{0.49\textwidth}
        \includegraphics[width=\textwidth]{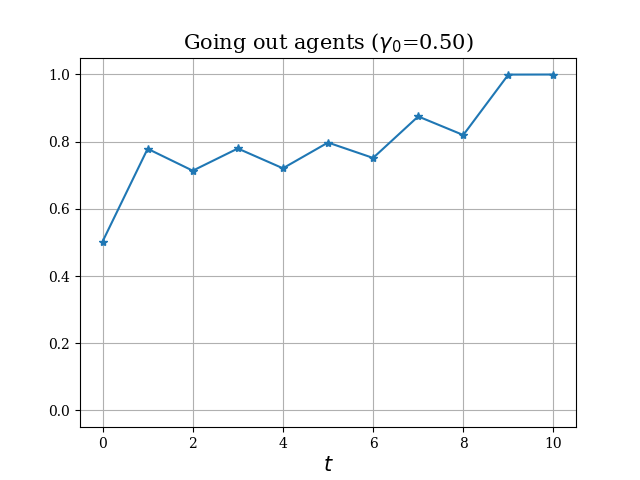}
    \end{minipage}
    \caption{Mean-field flow induced by the optimized policy under the population-level constraint \eqref{eq: macro_exp_con} with $\gamma_0 = 0.5$. The $x$-axis represents the time steps and the $y$-axis represents the percentage of infected agents (Left) or agents who go out (Right).}
    \label{fig: pop_ne}
\end{figure}

By Proposition~\ref{prop: pop cmfomo ne}, under population-level constraints, if a CMFG NE exists, then it must also be an unconstrained MFG NE. In Figure~\ref{fig: pop_ne}, we plot the mean-field flow induced by the optimized policy when $\gamma_0 = 0.5$ which is large enough such that a CMFG NE exists. By observation, the trajectories in Figure~\ref{fig: pop_ne} coincide with the trajectories corresponding to $\gamma_0 \geq 0.5$ in Figure~\ref{fig: SIS_sensitivity}. This is intuitive, since both the agent-population-level constraint \eqref{eq: gamma_0.25_con} and the population-level constraint \eqref{eq: macro_exp_con} become redundant when the value of $\gamma_0$ is large enough.

\bibliographystyle{unsrt}
\bibliography{refs.bib}

\appendix
\section{Proofs of results in \S\ref{sec: cmfomo}} \label{app: proofs}

\subsection*{Proof of Theorem~\ref{thm: suboptimal cmfomo}}

\begin{proof}
    We shall first prove the upper bound of the optimality gap. This part will be divided into three steps.

    \medskip
    \noindent
    \textbf{$G_{\text{opt}}$ step 1: bound the gap between $L$ and $\Psi(\pi)$.}
    
    By the same argument as the proof of \cite[Theorem 8]{mfomo}, we have:
    \begin{equation*}
        ||\Psi_t(\pi) - L_t||_1 \leq \frac{\epsilon}{c_2} \frac{(C_p + 1)^{t + 1} - 1}{C_p} \sqrt{|\mathcal{S}|}, \quad \forall t \in \T, 
    \end{equation*}
    which implies:
    \begin{equation*}
        ||\Psi(\pi) - L||_1 \leq \frac{\epsilon}{c_2} \frac{(C_p + 1)^{|\T| + 1} - 2C_p - 1}{C_p^2} \sqrt{|\mathcal{S}|} =: \epsilon \alpha(c_2).
    \end{equation*}

    \medskip
    \noindent
    \textbf{$G_{\text{opt}}$ step 2: near-optimality of $(\Psi(\pi), y, z, \lambda, w)$.}
    
    After replacing $L$ with $\Psi(\pi)$ in CMFOMO's objective \eqref{eq: cmfomo obj}, we have the following estimates. For notational simplicity, we write $\Psi$ in the place of $\Psi(\pi)$.
    \begin{align*}
        \left|\left|A^\top_\Psi y + z + r_\Psi - c_\Psi^\top \lambda \right|\right|_2 & \leq \left|\left|A^\top_\Psi y + z + r_\Psi - c_\Psi^\top \lambda \right|\right|_1 \\
        & \leq \left|\left|(A_\Psi - A_L)^\top y\right|\right|_1 + ||r_\Psi - r_L||_1 + \left|\left|(c_\Psi - c_L)^\top \lambda\right|\right|_1 + \frac{\epsilon}{c_1} \\
        & \leq ||\Psi - L||_1 \big(C_p ||y||_1 + C_r  + C_c ||\lambda||_1 \big) + \frac{\epsilon}{c_1} \\ 
        & \leq \epsilon \left[\frac{1}{c_1} + \alpha(c_2) \left( C_p \beta_y + C_r + k C_c \frac{|\T| r_{\max}}{\delta} \right)\right] =: \epsilon \zeta_1(c_1, c_2), \\
        A_\Psi \Psi - b & = 0, \\
        0 \leq z^\top \Psi & \leq \epsilon \left[\frac{1}{c_3} + \beta_z \alpha(c_2) \right] =: \epsilon \zeta_2(c_2, c_3), \\
        ||\gamma_0 - c_\Psi \Psi - w||_2 & \leq ||c_\Psi \Psi - c_L L||_2 + \frac{\epsilon}{c_4} \\
        & \leq ||c_\Psi (\Psi - L)||_2 + ||(c_\Psi - c_L) L||_2 + \frac{\epsilon}{c_4} \\
        & \leq \epsilon \left[\frac{1}{c_4} + \alpha(c_2) \Big(k c_{\max}  + C_c |\mathcal{S}| |\mathcal{A}| |\T|\Big)\right] =: \epsilon \zeta_3(c_2, c_4), \\
        \left|\lambda^\top (\gamma_0 - c_\Psi \Psi)\right| & \leq \left|\lambda^\top (c_\Psi \Psi - c_L L)\right| + \frac{\epsilon}{c_5} \\
        & \leq \epsilon \left[ \frac{1}{c_5} + k \alpha(c_2) \frac{|\T| r_{\max}}{\delta} \Big(k c_{\max} + C_c |\mathcal{S}| |\mathcal{A}| |\T|\Big) \right] \\ 
        & =: \epsilon \zeta_4(c_2, c_5).
    \end{align*}

    \medskip
    \noindent
    \textbf{$G_{\text{opt}}$ step 3: estimate the upper bound.}
    
    Throughout this step of the proof, we shall always take the mean-field flow $\Psi(\pi)$. Therefore, we omit the subscript for notational simplicity. 
    
    To start with, define $\Delta := A^\top y + z + r - c^\top \lambda$. Consider the following maximization problem:
    \begin{align*}
        \text{maximize}_{\hat{y}, \hat{z}, \hat\lambda} & \quad b^\top \hat{y} - \gamma_0^\top \hat\lambda \\
        \text{subject to} & \quad A^\top \hat{y} + \hat{z} = c^\top \hat\lambda - r + \Delta, \quad \hat\lambda \geq 0, \quad \hat{z} \geq 0,
    \end{align*}
    whose optimal solution is denoted by $\left(\hat{y}^\ast, \hat{z}^\ast, \hat{\lambda}^\ast\right)$. The dual problem is:
    \begin{align*}
        \text{minimize}_{\hat{d}} & \quad (-r + \Delta)^\top \hat{d} \\
        \text{subject to} & \quad A \hat{d} = b, \quad c \hat{d} \leq \gamma_0, \quad \hat{d} \geq 0,
    \end{align*}
    whose optimal solution is denoted by $\hat{d}^\ast$. In addition, denote by $d^\ast$ the optimal solution to the LP problem defined in Lemma~\ref{lemma: lp} driven by $\Psi(\pi)$. Then, we have the estimates:
    \begin{align}
        G_{\text{opt}}(\pi) := r^\top d^\ast - r^\top \Psi(\pi) & \leq r^\top d^\ast + \left(b^\top \hat{y}^\ast - \gamma_0^\top \hat{\lambda}^\ast\right) - \left(b^\top y - \gamma_0^\top \lambda\right) - r^\top \Psi(\pi) \nonumber \\
        & \leq r^\top d^\ast - (r - \Delta)^\top \hat{d}^\ast - \left[\left(b^\top y - \gamma_0^\top \lambda\right) + r^\top \Psi(\pi) \right] \nonumber \\
        & \leq ||\Delta||_1 + \epsilon \Big(|\T| \zeta_1 + \zeta_2 + \zeta_4 \Big) \label{eq: proof epxl upper bound} \\
        & \leq \epsilon \Big[(|\T| + 1) \zeta_1 + \zeta_2 + \zeta_4 \Big]. \nonumber
    \end{align}
    where the inequality \eqref{eq: proof epxl upper bound} comes from the estimate:
    \begin{align*}
        \left|\left(b^\top y - \gamma_0^\top \lambda\right) + r^\top \Psi(\pi)\right| & = \left|\Psi(\pi)^\top A^\top y - \gamma_0^\top \lambda + r^\top \Psi(\pi)\right| \\
        & \leq \left|\Psi(\pi)^\top (A^\top y + z + r - c^\top \lambda)\right| + \Psi(\pi)^\top z + \left|\lambda^\top (\gamma_0 - c \Psi(\pi))\right| \\
        & \leq \epsilon \Big(|\T| \zeta_1 + \zeta_2 + \zeta_4 \Big).
    \end{align*}
    This finishes our proof for the upper bound of the optimality gap.
    
    Next, we shall prove the upper bound of the feasibility gap. 
    
    \medskip
    \noindent
    \textbf{$G_{\text{fea}}$: estimate the upper bound.}
    
    Again, for notational simplicity, we write $\Psi$ in the place of $\Psi(\pi)$. It is enough to notice that:
    \begin{align*}
        G_{\text{fea}}(\pi) & = ||\min(0, \gamma_0 - c_\Psi \Psi)||_2 \\
        & \leq ||\min(0, \gamma_0 - c_\Psi \Psi) - \min(0, \gamma_0 - c_L L)||_2 + ||\min(0, \gamma_0 - c_L L)||_2 \\
        & \leq ||c_\Psi \Psi - c_L L||_2 + ||\min(0, \gamma_0 - c_L L)||_2 \\
        & \leq ||c_\Psi \Psi - c_L L||_2 + ||\gamma_0 - c_L L - w||_2 \\
        & \leq \epsilon \zeta_3. 
    \end{align*}
    This finishes our proof for the upper bound of the feasibility gap.

    Finally, we shall prove the lower bound of the optimality gap.
    
    \medskip
    \noindent
    \textbf{$G_{\text{opt}}$: estimate the lower bound.}

    In the rest, we shall always take the mean-field flow $\Psi(\pi)$. Therefore, we omit the subscript for notational simplicity. 

    To start with, define $\tilde\Delta := \max\left(0, c^\top \Psi(\pi) - \gamma_0\right)$. Consider the following maximization problem:
    \begin{align*}
        \text{maximize}_{\tilde{y}, \tilde{z}, \tilde\lambda} & \quad b^\top \tilde{y} - \left(\gamma_0 + \tilde\Delta\right)^\top \tilde\lambda \\
        \text{subject to} & \quad A^\top \tilde{y} + \tilde{z} = c^\top \tilde\lambda - r, \quad \tilde\lambda \geq 0, \quad \tilde{z} \geq 0,
    \end{align*}
    whose optimal solution is denoted by $\left(\tilde{y}^\ast, \tilde{z}^\ast, \tilde{\lambda}^\ast\right)$. The dual problem is:
    \begin{align*}
        \text{minimize}_{\tilde{d}} & \quad -r^\top \tilde{d} \\
        \text{subject to} & \quad A \tilde{d} = b, \quad c \tilde{d} \leq \gamma_0 + \tilde\Delta, \quad \tilde{d} \geq 0,
    \end{align*}
    whose optimal solution is denoted by $\tilde{d}^\ast$. In addition, denote by $\left(d^\ast, y^\ast, z^\ast, \lambda^\ast\right)$ a solution of $\K(\Psi(\pi))$ (\textit{i.e.}, the KKT conditions induced by $\Psi(\pi)$). Then, we have the estimates:
    \begin{align*}
        G_{\text{opt}}(\pi) := r^\top d^\ast - r^\top \Psi(\pi) & \geq r^\top d^\ast - r^\top \tilde{d}^\ast \\
        & = b^\top \tilde{y}^\ast - \left(\gamma_0 + \tilde\Delta\right)^\top \tilde{\lambda}^\ast - \left(b^\top y^\ast - \gamma_0 \lambda^\ast \right) \\
        & \geq b^\top y^\ast - \left(\gamma_0 + \tilde\Delta\right)^\top \lambda^\ast - \left(b^\top y^\ast - \gamma_0 \lambda^\ast \right) \\
        & = -\tilde{\Delta}^\top \lambda^\ast \\
        & \geq - \epsilon \sqrt{k} \frac{|\T| \rmax}{\delta} \zeta_3,
    \end{align*}
    where the last inequality comes from the facts that $G_{\text{fea}}(\pi) \leq \epsilon \zeta_3$ and $\left|\left|\lambda^\ast\right|\right|_\infty \leq \frac{|\T| \rmax}{\delta}$ (see Lemma~\ref{lemma: dual selection}). This finalizes our proof.
\end{proof}

\subsection*{Proof of Theorem~\ref{thm: gamma_eps}}

\begin{proof}
    The feasibility gap is a direct corollary of Theorem~\ref{thm: suboptimal cmfomo}, which further implies that the lower bound of the optimality gap is 0. Next, we prove the upper bound of the optimality gap. The proof utilizes similar techniques as the proof of Theorem~\ref{thm: suboptimal cmfomo}.
    
    First, the inequalities in $\opt$ step 2 of the proof of Theorem~\ref{thm: suboptimal cmfomo} become:
    \begin{align*}
        \left|\left|A^\top_\Psi y + z + r_\Psi - c_\Psi^\top \lambda \right|\right|_1 & \leq \frac{\epsilon_0}{\zeta_3} \left[\frac{1}{c_1} + \alpha(c_2) \left( C_p \tilde{\beta}_y + C_r + k C_c \frac{|\T| r_{\max}}{\delta - \epsilon_0} \right)\right] =: \frac{\epsilon_0}{\zeta_3} \tilde{\zeta}_1(c_1, c_2), \\
        A_\Psi \Psi - b & = 0, \\
        0 \leq z^\top \Psi & \leq \frac{\epsilon_0}{\zeta_3} \left[\frac{1}{c_3} + \tilde{\beta}_z \alpha(c_2) \right] =: \frac{\epsilon_0}{\zeta_3} \tilde{\zeta}_2(c_2, c_3), \\
        ||\gamma_0 - c_\Psi \Psi - w||_2 & \leq \frac{\epsilon_0}{\zeta_3} \left[\frac{1}{c_4} + \alpha(c_2) \Big(k c_{\max}  + C_c |\mathcal{S}| |\mathcal{A}| |\T|\Big)\right] + \epsilon_0 = \frac{\epsilon_0}{\zeta_3} \zeta_3(c_2, c_4) + \epsilon_0, \\
        \left|\lambda^\top (\gamma_0 - c_\Psi \Psi)\right| & \leq \frac{\epsilon_0}{\zeta_3} \left[ \frac{1}{c_5} + k \alpha(c_2) \frac{|\T| r_{\max}}{\delta - \epsilon_0} \Big(k c_{\max} + C_c |\mathcal{S}| |\mathcal{A}| |\T|\Big) \right] + \epsilon_0 \frac{k|\T|r_{\max}}{\delta - \epsilon_0} \\ 
        & =: \frac{\epsilon_0}{\zeta_3} \tilde{\zeta}_4(c_2, c_5) + \epsilon_0 \frac{k|\T|r_{\max}}{\delta - \epsilon_0}.
    \end{align*}
    Then, the upper bound of the optimality gap is obtained by repeating exactly the same calculation as $\opt$ step 3 (plugging the upper bounds of the above inequalities). The whole proof completes.
\end{proof}

\subsection*{Proof of Theorem~\ref{thm: char}}

\begin{proof}
    We shall always take the mean-field flow $L = \Psi(\pi)$. Therefore, we omit the subscript for notational simplicity. 

    To start with, notice that $A L = b$, which comes from the fact that $L = \Psi(\pi)$. Then, taking advantage of Lemma~\ref{lemma: dual selection}, direct calculations give the following estimates:
    \begin{align*}
        0 \leq z^\top L & = (c^\top \lambda - A^\top y - r)^\top L \\
        & = \lambda^\top \gamma_0 +\lambda^\top (cL - \gamma_0) - y^\top b - V^\pi(L) \\
        & = \lambda^\top \gamma_0 + \lambda^\top (cL - \gamma_0) + \left(r_L - c^\top \lambda\right)^\top d - V^\pi(L) \\
        & = \lambda^\top (cL - \gamma_0) + V^\ast_c(L) - V^\pi_{\mu_0}(L) \\
        & \leq \lambda^\top \max(0, cL - \gamma_0) + |\opt(\pi)| \\
        & \leq \sqrt{k}\frac{|\T|\rmax}{\delta} \epsilon_2 + \epsilon_1, \\
        ||\gamma_0 - cL - w||_2 & = ||\min(0, \gamma_0 - cL)||_2 \leq \epsilon_2, \\
        \left|\lambda^\top (\gamma_0 - c L)\right| & = \left|\lambda^\top c (d - L)\right| \\
        & = \left|y^\top A (d - L) + z^\top (d - L) + r^\top (d - L)\right| \\
        & = \left|-z^\top L + r^\top (d - L)\right| \\
        & \leq z^\top L + |\opt(\pi)| \\
        & \leq \sqrt{k}\frac{|\T|\rmax}{\delta} \epsilon_2 + 2\epsilon_1.
    \end{align*}
    Putting everything together yields the desired result.
\end{proof}

\section{Proofs of results in \S\ref{sec: N_player}} \label{app: N_player}

With an abuse of notation, denote by $\cmdp{L, \gamma_0}$ the constrained MDP problem defined in \S\ref{sec: setup} with the driving mean-field flow $L$ and threshold $\gamma_0$ (\textit{i.e.}, we explicitly write out the threshold $\gamma_0$ in $\cmdp{L}$), and denote by $\optimal{L, \gamma_0}$ its maximum cumulative expected reward, if exists.

\begin{lemma} \label{lemma: sensitivity_gamma}
     Fix the mean-field flow $L$. Let $\epsilon > 0$ and $\gamma_1, \gamma_2 \in \R^k$ be such that $||\gamma_1 - \gamma_2||_1 \leq \epsilon$. If the constrained MDP problems $\cmdp{L, \gamma_1}$ and $\cmdp{L, \gamma_2}$ satisfy Assumption~\ref{ass: strict feasibility}, then their optimal values will satisfy:
    \begin{equation*}
        |\optimal{L, \gamma_1} - \optimal{L, \gamma_2}| \leq \epsilon \frac{|\T| r_{\max}}{\delta}.
    \end{equation*}
\end{lemma}

\begin{proof}
    By Lemma~\ref{lemma: lp}, the MDP problems are equivalent to the following LP problems ($i = 1, 2$):
    \begin{align*}
        \text{minimize}_{d} & \quad -r_L^\top d \\
        \text{subject to} & \quad A_L d = b, \quad c_L d \leq \gamma_i, \quad d \geq 0,
    \end{align*}
    whose dual problems are:
    \begin{align*}
        \text{maximize}_{y, z, \lambda} & \quad b^\top y - \gamma_i^\top \lambda \\
        \text{subject to} & \quad A_L^\top y + z - c_L^\top \lambda = - r_L, \quad z \geq 0, \quad \lambda \geq 0.
    \end{align*}
    Denote by $(d_i^\ast, y_i^\ast, z_i^\ast, \lambda_i^\ast)$ the primal and dual optimal solutions to the above LP problems for $i = 1, 2$. By strong duality, we have:
    \begin{align*}
        -r_L^\top d_1^\ast & = b^\top y_1^\ast - \gamma_1^\top \lambda_1^\ast \\
        & = b^\top y_1^\ast - \gamma_2^\top \lambda_1^\ast + (\gamma_2^\top \lambda_1^\ast - \gamma_1^\top \lambda_1^\ast) \\
        & \leq b^\top y_2^\ast - \gamma_2^\top \lambda_2^\ast + (\gamma_2^\top \lambda_1^\ast - \gamma_1^\top \lambda_1^\ast) \\
        & = -r_L^\top d_2^\ast + (\gamma_2^\top \lambda_1^\ast - \gamma_1^\top \lambda_1^\ast) \\
        & \leq -r_L^\top d_2^\ast + \epsilon \frac{|\T| r_{\max}}{\delta},
    \end{align*}
    where the last inequality comes from Lemma~\ref{lemma: dual selection}. By exactly the same steps, we have:
    \begin{equation*}
        -r_L^\top d_2^\ast \leq -r_L^\top d_1^\ast + \epsilon \frac{|\T| r_{\max}}{\delta}.
    \end{equation*}
    This completes the proof.
\end{proof}

In the following, with an abuse of notation, denote by $L^{(N)}\left(\pi_1, \pi^{\otimes(N-1)}\right)$ the empirical distribution of the $N$ players when the first player takes policy $\pi_1$ and the rest $N-1$ players take policy $\pi$ (\textit{i.e.}, $\pi^{\otimes(N-1)}$ denotes $N - 1$ copies of $\pi$). Note that $L^{(N)}\left(\pi_1, \pi^{\otimes(N-1)}\right)$ is a random variable taking values in $[\Delta(\state \times \action)]^{|\T|}$. Recall that we embed $L^{(N)}$ in $\R^{|\state||\action||\T|}_{\geq 0}$, and for each $t \in \T$, we embed $L^{(N)}_t$ in $\R^{|\state||\action|}_{\geq 0}$.

\begin{lemma} \label{lemma: independence}
    Under Assumption~\ref{ass: continuous}, let $\pi$ and $\pi_1$ be two policies. We have:
    \begin{equation*}
        \max_{t \in \T} \E \left|\left| L^{(N)}_t\left(\pi_1, \pi^{\otimes(N-1)}\right) - \Psi_t(\pi) \right|\right|_1 \leq \left(\frac{1}{2\sqrt{N}} + \frac{2}{N}\right) |\state| |\action| \frac{C_{p, \state, \action}^{|\T|} - 1}{C_{p, \state, \action} - 1},
    \end{equation*}
    where the constant
    \begin{equation*}
        C_{p, \state, \action} := \left(C_p + 1\right) |\state| |\action|.
    \end{equation*}
\end{lemma}

\begin{proof}
    For notational simplicity, we use $L^{(N)}$ to represent $L^{(N)}_t\left(\pi_1, \pi^{\otimes(N-1)}\right)$ in this proof. Observe that for each $1 \leq t \leq T$ and $(s, a) \in \state \times \action$, we have:
    \begin{align*}
        \E \left|L^{(N)}_t(s, a) - \Psi_t(\pi)(s, a) \right| & = 
        \E \left[\E \left(\left.\left|L^{(N)}_t(s, a) - \Psi_t(\pi)(s, a) \right| \,\right|\, L^{(N)}_{t - 1} \right) \right] \\
        & = \E \left[\E \left(\left.\left| L^{(N)}_t(s, a) - \Big[\Psi_{t - 1}(\pi) \cdot p_{t - 1}\left(s \,|\,\Psi_{t-1}(\pi) \right) \Big] \pi_t (s)(a) \right| \,\right|\, L^{(N)}_{t - 1} \right) \right] \\
        & \leq \E(J_1 + J_2),
    \end{align*}
    where
    \begin{align*}
        J_1 & = \E \left(\left.\left| \frac{1}{N} \sum_{i = 1}^N \delta_{\left(s^{(i)}_t, a^{(i)}_t \right) = (s, a)} - \E\left( \frac{1}{N} \sum_{i = 1}^N \delta_{\left(s^{(i)}_t, a^{(i)}_t \right) = (s, a)} \right) \right| \,\right|\, L^{(N)}_{t - 1}\right) \leq \frac{1}{2\sqrt{N}}, \\
        J_2 & = \E \left(\left.\left| \E\left( \frac{1}{N} \sum_{i = 1}^N \delta_{\left(s^{(i)}_t, a^{(i)}_t \right) = (s, a)} \right) - \Big[\Psi_{t - 1}(\pi) \cdot p_{t - 1}\left(s \,|\,\Psi_{t-1}(\pi) \right)\Big] \pi_t (s)(a) \right| \,\right|\, L^{(N)}_{t - 1}\right) \\
        & \leq \frac{1}{N} \left| p_{t-1}\left(s^{(1)}_{t-1}, a^{(1)}_{t-1}, L^{(N)}_{t-1}\right)(s) \pi_{1, t}(s)(a) - \Big[\Psi_{t - 1}(\pi) \cdot p_{t - 1}\left(s \,|\,\Psi_{t-1}(\pi) \right)\Big] \pi_t (s)(a) \right| \nonumber \\
        & \quad + \frac{1}{N} \left| \sum_{i=2}^N \left[p_{t-1}\left(s^{(i)}_{t-1}, a^{(i)}_{t-1}, L^{(N)}_{t-1}\right)(s) \pi_t(s)(a) - \Big[\Psi_{t - 1}(\pi) \cdot p_{t - 1}\left(s \,|\,\Psi_{t-1}(\pi) \right)\Big] \pi_t (s)(a)\right] \right| \\
        & \leq \frac{1}{N} + \frac{1}{N} \left| \sum_{i=2}^N \left[ p_{t-1}\left(s^{(i)}_{t-1}, a^{(i)}_{t-1}, L^{(N)}_{t-1}\right)(s) - p_{t-1}\left(s^{(i)}_{t-1}, a^{(i)}_{t-1}, \Psi_{t-1}(\pi)\right)(s) \right] \right| \nonumber \\
        & \quad + \frac{1}{N} \left| \sum_{i=2}^N \left[ p_{t-1}\left(s^{(i)}_{t-1}, a^{(i)}_{t-1}, \Psi_{t-1}(\pi)\right)(s) - \Psi_{t - 1}(\pi) \cdot p_{t - 1}\left(s \,|\,\Psi_{t-1}(\pi) \right) \right] \right| \\
        & \leq \frac{1}{N} + \frac{N - 1}{N} C_p \left|\left| L_{t-1}^{(N)} - \Psi_{t-1}(\pi) \right|\right|_1 + \left(\left|\left| L_{t-1}^{(N)} - \Psi_{t-1}(\pi) \right|\right|_1 + \frac{1}{N} \right) \\
        & \leq \frac{2}{N} + (C_p + 1) \left|\left| L_{t-1}^{(N)} - \Psi_{t-1}(\pi) \right|\right|_1.
    \end{align*}
    When $t = 0$, by a similar argument as above, we get:
    \begin{align*}
        \E \left| L_0^{(N)}(s, a) - \Psi_0(\pi)(s, a) \right| & = \E \left| L_0^{(N)}(s, a) - \mu_0(s) \pi(s)(a) \right| \\
        & \leq \frac{1}{2\sqrt{N}} + \frac{1}{N}.
    \end{align*}
    The desired result is obtained by an induction argument on time stamps.
\end{proof}

Similarly, we denote by $L^{(N)}\left(\pi_1 \,\left|\, \pi^{\otimes (N-1)}\right.\right)$ the empirical distribution of the first player when the first player takes policy $\pi_1$ and the rest players take policy $\pi$. Note that $L^{(N)}\left(\pi_1 \,\left|\, \pi^{\otimes (N-1)}\right.\right)$ is a random variable taking values in $[\Delta(\state \times \action)]^{|\T|}$. Different from our previously defined $L^{(N)}\left(\pi_1, \pi^{\otimes(N-1)}\right)$, here $L^{(N)}\left(\pi_1 \,\left|\, \pi^{\otimes (N-1)}\right.\right)$ only characterizes the distribution of the first player. In addition, in the CMFG, we denote by $L(\pi_1 \,|\, \pi)$ the distribution of a specific player who takes policy $\pi_1$, provided that that population distribution is $\Psi(\pi)$ (\textit{i.e.}, all the rest players take policy $\pi$).

\begin{lemma} \label{lemma: first_agent_policy}
    Under Assumption~\ref{ass: continuous}, let $\pi$ and $\pi_1$ be two policies. Then, we have:
    \begin{equation*}
        \max_{t \in \T} \left|\left| \E \left[L^{(N)}_t\left(\pi_1 \,\left|\, \pi^{\otimes (N-1)}\right.\right)\right] - L_t(\pi_1 \,|\, \pi)\right|\right|_1 \leq \left(\frac{1}{2\sqrt{N}} + \frac{2}{N}\right) C_p |\state| |\action| (|\T| - 1)\frac{C_{p, \state, \action}^{|\T|} - 1}{C_{p, \state, \action} - 1},
    \end{equation*}
    where $C_{p, \state, \action}$ is defined in Lemma~\ref{lemma: independence}.
\end{lemma}

\begin{proof}
    For notational simplicity, in the rest of the proof, we use $L^{(N)}$ to represent $L^{(N)}\left(\pi_1 \,\left|\, \pi^{\otimes (N-1)}\right.\right)$, and use $L$ to represent $L(\pi_1 \,|\, \pi)$.
    
    For each $1 \leq t \leq T$ and $(s, a) \in \state \times \action$, we have:
    \begin{align*}
        & \left| \E \left[L^{(N)}_t\right] (s, a) - L_t(s, a) \right| \\
        & = \left| \E\left[ L_{t - 1}^{(N)} \cdot p_{t - 1} \left(s \,\left|\, L^{(N)}_{t-1}(\pi_1, \pi) \right.\right) \right] \pi_1(s)(a) - \Big[L_{t-1} \cdot p_{t-1} (s \,|\, \Psi_{t-1}(\pi))\Big] \pi_1(s)(a) \right| \\
        & \leq \left| \E\left[ L_{t - 1}^{(N)} \cdot p_{t - 1} \left(s \,\left|\, L^{(N)}_{t-1}(\pi_1, \pi) \right.\right) \right] - \Big[L_{t-1} \cdot p_{t-1} (s \,|\, \Psi_{t-1}(\pi))\Big] \right| \\
        & \leq \left| \E\left[ L_{t - 1}^{(N)} \cdot p_{t - 1} \left(s \,\left|\, L^{(N)}_{t-1}(\pi_1, \pi) \right.\right) - L_{t - 1}^{(N)} \cdot p_{t - 1} (s \,|\, \Psi_{t-1}(\pi)) \right]\right| \nonumber \\
        & \quad + \left| \E\left[L_{t - 1}^{(N)} \cdot p_{t - 1} (s \,|\, \Psi_{t-1}(\pi)) \right] - L_{t-1} \cdot p_{t-1} (s \,|\, \Psi_{t-1}(\pi)) \right| \\
        & \leq \left(\frac{1}{2\sqrt{N}} + \frac{2}{N}\right) C_p |\state| |\action| \frac{C_{p, \state, \action}^{|\T|} - 1}{C_{p, \state, \action} - 1} + \left|\left| \E \left[L^{(N)}_{t-1}\right] - L_{t-1}\right|\right|_1.
    \end{align*}
    We mention that the first term of the last inequality comes from Lemma~\ref{lemma: independence}, where $C_{p, \state, \action}$ is defined.
    
    Notice that when $t = 0$, we have:
    \begin{align*}
        \left|\left| \E \left[L^{(N)}_0\right] - L_0\right|\right|_1 = 0.
    \end{align*}
    The desired result is obtained by an induction argument on time stamps.
\end{proof}

With an abuse of notation, denote by $V^{\pi_1}\left(\pi^{\otimes(N-1)}\right)$ the expected cumulative reward of the first player who takes policy $\pi_1$, provided that the rest players take policy $\pi$, \textit{i.e.},
\begin{equation*}
    V^{\pi_1}\left(\pi^{\otimes(N-1)}\right) := \E_{a^{(1)} \sim \pi_1,\, \forall 2 \leq j \leq N,\, a^{(j)} \sim \pi} \left[ \sum_{t \in \T} r_t\left(s^{(1)}_t, a^{(1)}_t, L^{(N)}_t\right) \,\bigg|\, s^{(1)}_0 \sim \mu_0 \right].
\end{equation*}
On top of Lemma~\ref{lemma: independence} and Lemma~\ref{lemma: first_agent_policy}, a direct calculation produces the following estimate.

\begin{lemma} \label{lemma: reward_approx}
    Under Assumption~\ref{ass: continuous}, let $\pi$ and $\pi_1$ be two policies. Then, we have:
    \begin{equation*}
        \left| V^{\pi_1}\left(\pi^{\otimes(N-1)}\right) - V^{\pi_1}(\Psi(\pi)) \right| \leq \left(\frac{1}{2\sqrt{N}} + \frac{2}{N}\right) C_p (C_r + r_{\max}) |\state| |\action| (|\T| - 1)^2 \frac{C_{p, \state, \action}^{|\T|} - 1}{C_{p, \state, \action} - 1},
    \end{equation*}
    where $C_{p, \state, \action}$ is defined in Lemma~\ref{lemma: independence}.
\end{lemma}

Similarly, in the $N$-player game, we denote:
\begin{align*}
    \C^{\pi_1}\left(\pi^{\otimes(N-1)}\right) := \E_{a^{(1)} \sim \pi_1,\, \forall 2 \leq j \leq N,\, a^{(j)} \sim \pi} \left[ \sum_{t \in \T} c_t\left(s^{(1)}_t, a^{(1)}_t, L^{(N)}_t\right) \,\bigg|\, s^{(1)}_0 \sim \mu_0 \right],
\end{align*}
and in the CMFG, we denote (with an abuse of notation):
\begin{align*}
    \C^{\pi_1}\left(L\right) := \E_{a \sim \pi_1} \left[ \sum_{t \in \T} c_t\left(s_t, a_t, L_t\right) \,\bigg|\, s_0 \sim \mu_0 \right].
\end{align*}
By replacing $r$ with $c$ in Lemma~\ref{lemma: reward_approx}, we get the following estimate on the value of the cost functions.

\begin{lemma} \label{lemma: cost_approx}
    Under Assumption~\ref{ass: continuous}, let $\pi$ and $\pi_1$ be two policies. Then, we have:
    \begin{equation*}
        \left|\left| \C^{\pi_1}\left(\pi^{\otimes(N-1)}\right) - \C^{\pi_1}(\Psi(\pi)) \right|\right|_\infty \\
        \leq \left(\frac{1}{2\sqrt{N}} + \frac{2}{N}\right) C_p (C_c + c_{\max}) |\state| |\action| (|\T| - 1)^2 \frac{C_{p, \state, \action}^{|\T|} - 1}{C_{p, \state, \action} - 1},
    \end{equation*}
    where $C_{p, \state, \action}$ is defined in Lemma~\ref{lemma: independence}.
\end{lemma}

Finally, we are prepared to prove the main result Theorem~\ref{thm: eps-NE}.

\begin{proof}[Proof of Theorem~\ref{thm: eps-NE}]
    Fix $\epsilon$ and $N$ that satisfy the conditions. In the rest of the proof, we let all players in the $N$-player game take the policy $\pi^\ast$. Due to the symmetry of the players, it is enough to analyze the first player.
    
    First, by Lemma~\ref{lemma: cost_approx}, it is not hard to see that $\F^{(1)}$ is non-empty (\textit{e.g.}, the policy in Assumption~\ref{ass: strict feasibility N} belongs to $\F^{(1)}$). Then, we have:
    \begin{align*}
        \sup_{\pi_1 \in \mathcal{F}^{(1)}} V^{\pi_1}\left((\pi^\ast)^{\otimes(N-1)}\right) & = V^{\hat{\pi}_1}\left((\pi^\ast)^{\otimes(N-1)}\right) \leq V^{\hat{\pi}_1}(\Psi(\pi^\ast)) + \epsilon,
    \end{align*}
    where we denote that the supremum is achieved at $\hat{\pi}_1$. Since $\hat{\pi}_1 \in \mathcal{F}^{(1)}$, by Lemma~\ref{lemma: cost_approx}, we have:
    \begin{align*}
        \fea(\hat{\pi}_1, \Psi(\pi^\ast)) \leq \sqrt{k} \epsilon,
    \end{align*}
    where $\fea$ is defined in Definition~\ref{def: fea_gap}. Therefore, Lemma~\ref{lemma: sensitivity_gamma} gives the following estimate:
    \begin{align*}
        V^{\hat{\pi}_1}(\Psi(\pi^\ast)) \leq V^{\pi^\ast}(\Psi(\pi^\ast)) + \epsilon \frac{k |\T| r_{\max}}{\delta}.
    \end{align*}
    Putting everything together, we get:
    \begin{align*}
        \sup_{\pi_1 \in \mathcal{F}^{(1)}} V^{\pi_1}\left((\pi^\ast)^{\otimes(N-1)}\right) & \leq V^{\pi^\ast}(\Psi(\pi^\ast)) + \epsilon \left(1 + \frac{k |\T| r_{\max}}{\delta}\right) \\
        & \leq V^{\pi^\ast}\left((\pi^\ast)^{\otimes(N-1)}\right) + \epsilon \left(2 + \frac{k |\T| r_{\max}}{\delta}\right).
    \end{align*}
    This completes the proof of the $\epsilon_1$ part. The $\epsilon_2$ part is a direct corollary of Lemma~\ref{lemma: cost_approx}. This completes the whole proof.
\end{proof}

\end{document}

%% file: defs_cmfg.tex
\newcommand{\BEAS}{\begin{eqnarray*}}
\newcommand{\EEAS}{\end{eqnarray*}}
\newcommand{\BEQ}{\begin{equation}}
\newcommand{\EEQ}{\end{equation}}
\newcommand{\BIT}{\begin{itemize}}
\newcommand{\EIT}{\end{itemize}}

\newcommand{\E}{\mathbb{E}}

\def\<#1,#2>{\langle #1,#2\rangle}

\newtheorem{theorem}{Theorem}[section]
\newtheorem{remark}{Remark}[section]
\newtheorem{assumption}{Assumption}
\newtheorem{lemma}[theorem]{Lemma}
\newtheorem{corollary}[theorem]{Corollary}
\newtheorem{proposition}[theorem]{Proposition}
\newtheorem{definition}{Definition}
{\begin{quote}}{\end{quote}}

\makeatletter
\long\def\@makecaption#1#2{
   \vskip 9pt
   \begin{small}
   \setbox\@tempboxa\hbox{{\bf #1:} #2}
   \ifdim \wd\@tempboxa > 5.5in
        \begin{center}
        \begin{minipage}[t]{5.5in}
        \addtolength{\baselineskip}{-0.95pt}
        {\bf #1:} #2 \par
        \addtolength{\baselineskip}{0.95pt}
        \end{minipage}
        \end{center}
   \else
	\hbox to\hsize{\hfil\box\@tempboxa\hfil}
   \fi
   \end{small}\par
}
\makeatother

\newcounter{oursection}

\newcounter{lecture}

\newcommand{\action}{\mathcal{A}}
\newcommand{\C}{\mathcal{C}}
\newcommand{\F}{\mathcal{F}}
\newcommand{\K}{\mathcal{K}}

\newcommand{\R}{\mathbb{R}}
\newcommand{\state}{\mathcal{S}}
\newcommand{\T}{\mathcal{T}}

\newcommand{\rmax}{r_{\max}}
\newcommand{\cmax}{c_{\max}}
\newcommand{\cmdp}[1]{\mathcal{M}_c (#1)}
\newcommand{\mdp}[1]{\mathcal{M} (#1)}
\newcommand{\cmdpN}[1]{\mathcal{M}_c^{(#1)}}
\newcommand{\optimal}[1]{V^\ast_c (#1)}
\newcommand{\rv}{\mathbf}

\newcommand{\opt}{G_{\text{opt}}}
\newcommand{\fea}{G_{\text{fea}}}

\theoremstyle{definition}
\newtheorem*{cmdps}{Constrained MDP (CMDP)}
\newtheorem*{mdps}{Unconstrained MDP}
\newtheorem*{N_cmdps}{Constrained Control Problem of the $i$-th Player}

\newtheorem*{kkt}{KKT Conditions}
\newtheorem*{popkkt}{KKT Conditions (Population-level)}
\newtheorem{example}{Example}